\documentclass{article}
\usepackage[a4paper, total={5.5in, 8in}]{geometry}
\usepackage{amsfonts, amsmath, amsthm, dsfont}
\usepackage[utf8x,utf8]{inputenc}
\usepackage{centernot}
\usepackage{color}
\usepackage{comment}
\usepackage{tikz}
\usepackage[shortlabels]{enumitem}  
\usetikzlibrary{shapes.geometric, arrows}
\usepackage{empheq}

\usepackage[colorlinks=true,urlcolor=black]{hyperref}
\definecolor{refcol}{rgb}{0.1,0,0.6}
\hypersetup{colorlinks}
\hypersetup{linkcolor=refcol, citecolor=refcol}

\usepackage{algorithm}
\usepackage[noend]{algpseudocode}
\usepackage{subcaption}

\renewcommand{\d}{\textnormal{d}}

\newcommand{\E}{\mathbb{E}}
\newcommand{\Indicator}[1]{\mathds{1}\left(#1\right)}
\renewcommand{\P}{\mathbb{P}}
\newcommand{\R}{\mathbb{R}}

\newcommand{\Var}{\mathbb{V}\text{\textnormal{ar}}}

\newtheorem{lemma}{Lemma}[section]
\newtheorem{corollary}[lemma]{Corollary}
\newtheorem{thm}[lemma]{Theorem}
\newtheorem{rk}[lemma]{Remark}
\newtheorem{prop}[lemma]{Proposition}
\newtheorem{defn}[lemma]{Definition}

\newtheorem{assumption}[lemma]{Assumption}

\newcommand{\indicatorA}[1]{ \mathds{1}_{A}(#1) }

\usepackage[textsize=tiny]{todonotes}
\definecolor{jb}{rgb}{0.8,0.1,0.1}

\definecolor{sw}{rgb}{0.1,0.5,0.1}

\newcommand{\eps}{\varepsilon}

\newcommand{\N}{{\mathbb N}}

\newcommand{\bX}{{\mathbf X}}

\newcommand{\grad}{\nabla}

\newcommand{\rhosig}{\rho_{\sigma}^{*}}
\newcommand{\rhoopt}{\rho^{*}}

\newcommand{\Hess}{\textnormal{Hess}}
\newcommand{\vMFNM}{\textnormal{vMFNM}}
\newcommand{\vMFN}{\textnormal{vMFN}}
\newcommand{\vMF}{\textnormal{vMF}}
\newcommand{\cumu}{\textnormal{cumu}}

\usepackage{titlesec}

\titleformat{\subsection}[runin]
       {\normalfont\bfseries}
       {\thesubsection}
       {0.5em}
       {}
       [.]
\usepackage[style=alphabetic,doi=false,isbn=false,url=false,eprint=true,giveninits=true ,maxcitenames=9,maxbibnames=9]{biblatex}
\AtEveryBibitem{%
  \clearlist{language}%
  \clearfield{month}
  \clearfield{eprint}
}
\title{Affine invariant interacting Langevin dynamics in Markov chain importance sampling for rare event estimation}

\usepackage{authblk}
\author[1]{Jason Beh}
\author[1]{Jérôme Morio}
\author[2]{Florian Simatos}
\author[3]{Simon Weissmann}
\affil[1]{ONERA/DTIS, Fédération ENAC ISAE-SUPAERO ONERA, Université de Toulouse, F-31055 Toulouse, France}
\affil[2]{Fédération ENAC ISAE-SUPAERO ONERA, Université de Toulouse, Toulouse, France}
\affil[3]{Universität Mannheim, Mathematisches Institut, D-68159 Mannheim, Germany}
\date{}

\begin{document}

\maketitle
\begin{abstract}
    This work considers the framework of Markov chain importance sampling~(MCIS), 
    in which one employs a Markov chain Monte Carlo~(MCMC) scheme to sample particles approaching the optimal distribution for importance sampling, prior to estimating the quantity of interest through importance sampling. In rare event estimation, the optimal distribution admits a non-differentiable log-density, thus gradient-based MCMC can only target a smooth approximation of the optimal density. We propose a new gradient-based MCIS scheme for rare event estimation, called affine invariant interacting Langevin dynamics for importance sampling~(ALDI-IS), in which the affine invariant interacting Langevin dynamics~(ALDI) is used to sample particles according to the smoothed zero-variance density. We establish a non-asymptotic error bound when importance sampling is used in conjunction with samples independently and identically distributed according to the smoothed optiaml density to estimate a rare event probability, and an error bound on the sampling bias when a simplified version of ALDI, the unadjusted Langevin algorithm, is used to sample from the smoothed optimal density.
    We show that the smoothing parameter of the optimal density has a strong influence and exhibits a trade-off between a low importance sampling error and the ease of sampling using ALDI. We perform a numerical study of ALDI-IS and illustrate this trade-off phenomenon on standard rare event estimation test cases.
\end{abstract}
\tableofcontents
\titleformat{\section}[runin]
       {\normalfont\bfseries}
       {\thesection}
       {0.5em}
       {}
       [.]
\section{Introduction}
Consider a probability triple $(\Omega, \mathcal{F}, \P)$. In the context of static rare event estimation in $\R^d$~\cite{morio_survey_2014,cerou_adaptive_2019}, equipped with the Borel sigma-algebra $\mathcal{B}(\R^d)$, one is concerned with the estimation of the probability $p = \P(g(Z) \leq 0)$ where $Z: \Omega \to \R^d$ is distributed according to a distribution $\rho$ on $\R^d$, and $g: \R^d \to \R$ is $\mathcal{F}$-measurable. In the applicative context, $g(Z)$ is seen as the output of a blackbox system taking $d$-dimensional random input $Z$ with known distribution $\rho$, and one is interested in estimating~$\rho(A) = \P(g(Z) \leq 0)$, the measure of the failure set $A =\{z \in \R^d: g(z) \leq 0\}\in\mathcal B(\R^d)$ under~$\rho$. The function $g$ is then termed the limit state function~(LSF). Often, the gradient of $g$, denoted as $\nabla g$, is accessible in the output of the system. Owing to the practical knowledge on $\rho$, isoprobabilistic transformations $\mathcal{T}$, such as the Nataf or the Rosenblatt transformation~\cite{lebrun09b,lebrun09c}, can be applied on $Z$ so that $\mathcal{T}(Z)$ is normally distributed, and one views instead $g \circ \mathcal{T}^{-1}$ as the measurable function defining the failure set. As such, the design of a static rare event estimation algorithm in the standard normal space is common practice, and without loss of generality we hereby let $\rho$ be the $d$-dimensional standard normal distribution $N(0,I)$. For any measure absolutely continuous with respect to~(w.r.t.)\ the Lebesgue measure on $\R^d$, we deliberately use the same notation for its Radon-Nikodym derivative w.r.t.\ to the Lebesgue measure on $\R^d$, which we will exchangeably call its probability density function~(PDF).

\subsection{Importance sampling in rare event estimation} \label{subs:Intro-IS}
We denote the indicator function of a set $B \in \mathcal{B}(\R^d)$ as $\mathds{1}_B$, where for any $x \in \R^d$,
\[\mathds{1}_B(x) = \Indicator{x \in B}= \begin{cases}
    1 \ \text{ if } x\in B,
    \\0 \ \text{ otherwise.}
\end{cases}\]
The crude Monte Carlo estimator considers $N$ samples $(Z_i)_{i=1,\ldots,N}$ independently and identically distributed~(i.i.d.) according to $\rho$, and estimates $p$ as
\[\hat p_{\rho} = \frac{1}{N}\sum_{i=1}^N \mathds{1}_{A}(Z_i), \text{ with } \E[\hat p_{\rho}] = p \text{ and } \Var[\hat p_{\rho}] = \frac{p(1-p)}{N}. \]
Although it is unbiased, the normalized root mean square error~(nRMSE) of the crude Monte Carlo estimator, which coincides with its coefficient of variation because of the unbiasedness, writes as
\begin{align} \label{eq:nRMSE-MC}
 \text{nRMSE}(\hat p_{\rho}) := \frac{1}{p}\E\left[(\hat p_{\rho} - p )^2\right]^{1/2} = \left(\frac{1-p}{N p}\right)^{1/2}.
\end{align}
One then needs roughly $\varepsilon^{-2} p^{-1}$ samples to achieve an nRMSE of order $\varepsilon$, which amounts to an unreasonable amount of samples if the probability to estimate $p$ is low, in practice below $10^{-3}$.

The importance sampling~(IS) estimator instead considers an auxiliary distribution $\rho'$, w.r.t.\ which $\mathds{1}_A \rho$ is absolutely continuous. With $N$ samples $(X_i)_{i=1,\ldots,N}$ i.i.d.\ according to~$\rho'$, the estimator is defined as 
\begin{align}\label{def:IS-estimator}
\hat p_{\rho'} = \frac{1}{N}\sum_{i=1}^N \frac{\mathds{1}_A(X_i)\rho(X_i)}{\rho'(X_i)}, \text{ with } \E[\hat p_{\rho'}] = p.
\end{align}
The IS estimator 
remains unbiased, but its nRMSE depends crucially on the choice of $\rho'$. The zero-variance IS distribution $\rhoopt$, defined as 
\[ \rho^* = \frac{\mathds{1}_A \rho}{p},\]
makes the IS estimator $\hat p_{\rho^*}$ constant equal to $p$. 
We will refer to this distribution as the optimal distribution owing to its zero-variance property when used as the auxiliary distribution in importance sampling. However, the immediate use of the optimal distribution~$\rho^*$ is impossible, since its PDF involves the probability one seeks to estimate as its normalizing constant, and the importance sampling weights in the estimator~\eqref{def:IS-estimator} are $p$.

Topics in importance sampling for rare event estimation concern how to approach the optimal distribution $\rho^*$ in a practical manner. They can broadly be categorized in three ways. The first considers a parametric family of distributions, in which one finds the parameters which minimize a fixed metric w.r.t.\ the optimal distribution. The usual choice of metric is the Kullback-Leibler divergence since it defines a convex minimization problem if the parametric family is of the exponential-type family. The metric choice of the Kullback-Leibler divergence constitutes the basis of the popular cross-entropy algorithm~\cite{rubinstein_cross_2004} for rare event estimation, in which the optimal parameters are estimated sequentially. Once the optimal parameters are estimated to construct the parametric distribution, $N$ samples are generated to estimate $p$. This algorithm led to much follow-up works including numerical improvements~\cite{papaioannou_improved_2019,geyer_cross_2019,ELMASRI2021,uribe2021,el_masri2024,demange-chryst2024variational} and theoretical studies~\cite{homem-de-mello_study_2007,delyon21,beh2023insight}. 

The second category assumes specific structure on the rare event set $A$ and designs good auxiliary distributions based on some additional knowledge on $A$. Some examples in this category include highly reliable Markovian systems~\cite{LEcuyer2010,LEcuyer2011c}, static network reliability estimation~\cite{Cancela10,LEcuyer2011b}, large deviation estimation~\cite{LEcuyer2010,Guyader2020} and union of rare events~\cite{art19}. Our current work for static rare event estimation only considers the rare event set in the form of $A = \{ z \in \R^d: g(z) \leq 0 \}$ for which practitioners only have access to the evaluations of the LSF $g$ and its gradient $\grad g$. If in practice the rare event set $A$ satisfies the settings of some of the aforementioned works, their auxiliary densities designed specifically for $A$ should be employed~\cite{caron2012}.

The third category concerns the Markov chain importance sampling~(MCIS) framework, introduced in~\cite{Botev2013}. This framework generally approaches $\rho^*$ in the three following steps. 
\begin{enumerate}[label = {Step \arabic*.},align=left]
    \item In the Markov chain step, an MCMC algorithm with invariant distribution $\rho^*$ is first used to generate  $M$ samples from $\rho^*$. 
    \item To evaluate the PDF of the $M$ samples for constructing IS weights, a distribution is estimated from the MCMC samples using either a parametric, kernel-based, or nonparametric density estimator.
    \item In the final importance sampling step, $N$ i.i.d.~samples are generated from the fitted distribution, and an unbiased estimate of $p$ is obtained using the importance sampling estimator in~\eqref{def:IS-estimator}.
\end{enumerate}
In the present work, the focus is placed on the MCIS framework, for which we propose a new gradient-based scheme using affine-invariant interacting Langevin dynamics. Furthermore, we analyze the use of a smoothed optimal distribution for gradient-based MCIS, defined further in Section~\ref{subsection:smoothed-density}. 
\subsection{Literature overview on MCIS}
As mentioned, the MCIS scheme consists of three distinct steps, and we review related works for each step separately.

The general principle of MCMC is to construct a Markov chain to approach its invariant distribution, requiring only evaluations of the target density up to a normalizing constant. Consequently, MCMC algorithms are well-suited for sampling from $\rhoopt \propto \mathds{1}_A \rho$. The literature on MCMC is vast, for a recent survey, see, for instance, \cite{Luengo2020}. 

We highlight several MCMC methods that have been used in MCIS schemes for rare event estimation. In~\cite{Botev2013}, the Gibbs sampler~\cite{Gelfand01061990,robert_monte_2004} is employed, a special case of the Metropolis–Hastings algorithm that uses the full conditional distributions of the target as the transition kernel. In~\cite{schuster21}, although not focused on rare events, the authors use standard Metropolis–Hastings with a Gaussian random walk proposal, as well as the unadjusted Langevin algorithm (ULA), which we examine in more detail in Section~\ref{subsection:ULA}. The Metropolis-adjusted Langevin algorithm (MALA) has also been applied as MCIS scheme in~\cite{pmlr-v206-tit23a} and~\cite{cheng2025langevinbifidelityimportancesampling}. The former has a usual MCIS framework and the latter assumes access to a bi-fidelity code for the LSF, and MALA is run with low-fidelity code to save computational costs.

The authors in~\cite{PAPAKONSTANTINOU2023} explore Hamiltonian Monte Carlo (HMC), which requires evaluating the gradient of the log-density of the target distribution — defined as the logarithm of its PDF. A variant considered in the same work further requires the Hessian of the log-density. Since $\rhoopt$ includes the non-differentiable indicator function $\mathds{1}_A$, the authors introduce a smooth approximation of the indicator to enable gradient-based sampling. Such approximations are commonly employed when gradient information is needed in MCIS schemes.

Other algorithms that, while not classified as MCMC methods, still operate as particle-based samplers include several recent approaches. For example, \cite{wagner22} applies the ensemble Kalman filter by reformulating the task of sampling from $\rhoopt$ as a Bayesian inverse problem.
The ensemble Kalman filter was originally introduced as a particle-based approximation of the Bayesian filtering distribution in data assimilation \cite{Evensen2003}. More recently, it has been successfully applied to Bayesian inverse problems \cite{StLawIg2013}, which is the basis of~\cite{wagner22}.
Another example is provided in~\cite{althaus2024}, which considers consensus-based sampling (CBS), originally proposed in~\cite{Carrillo22}. CBS is inspired by the metaheuristic optimization technique known as consensus-based optimization~\cite{PTTM2017}. This sampling scheme defines a specific stochastic differential equation to govern particle dynamics, with the invariant distribution being a Laplace approximation (i.e., a Gaussian approximation) of the target distribution. Lastly,~\cite{ehre2024steinvariationalrareevent} applies the Stein Variational Gradient Descent (SVGD) algorithm~\cite{LW2016} to sample from the smoothed version $\rhoopt$. Like the Hamiltonian MCMC, SVGD requires gradient information of the target distribution.

Once the particles have been sampled in the MCMC step, their PDF must be approximated for use in the importance sampling estimation.
In rare event estimation, common parametric models include the Gaussian mixture model~\cite{wagner22,ehre2024steinvariationalrareevent} and the von Mises–Fisher–Nakagami mixture~(vMFNM) model~\cite{papaioannou_improved_2019,uribe2021}, the latter of which will be discussed in more detail in Section~\ref{subsection:vMFNM}. More recently, \cite{kruse2025scalableimportancesamplinghigh} proposed to fit a mixture model consisting of probabilistic principal component analyzers, which is a low-rank mixture model similarly to the vMFNM model. In~\cite{Botev2013}, a semi-parametric model based on the product of conditional densities is introduced for use with the Gibbs sampler. The use of nonparametric density estimators has also been investigated. In~\cite{Botev2013}, for instance, the Gibbs sampler's tractable conditional densities are exploited to construct a nonparametric estimate of the joint distribution. In the context of acceptance–rejection Metropolis–Hastings,~\cite{schuster21} suggest using all proposed particles, including rejected ones, to fit a nonparametric density based on the evaluated proposal densities, thereby avoiding further target density evaluations. In the context of general MCMC methods, \cite{zhang96} employs kernel density estimation to approximate the sample distribution. Meanwhile, \cite{chan_improved_2012} introduces an alternative to the standard expectation-maximization approach by fitting a parametric model whose parameters maximize the variance of the log-density, based on the insight that $\rhoopt$ as the zero-variance distribution, has an infinite log-density variance.
\subsection{Smoothed optimal density for gradient-based MCIS scheme} \label{subsection:smoothed-density}
Certain classes of MCMC methods such as the Hamiltonian MCMC require the evaluation of the gradient of the target density.
However, the target $\rho^*$ is non-differentiable on $\R^d$ due to the indicator function of the set $A$. To mitigate this issue, we consider the following smooth approximation of the indicator function. More precisely, we follow a smooth approximation as defined in~\cite{uribe2021,ehre2024steinvariationalrareevent} for the application to rare event estimation, where the indicator function is approximated based on the logistic function $x \mapsto 1/(1+e^{-x})$. We note that alternative choices of $F_\sigma$ can be found in~\cite{papaioannou_sequential_2016,papaioannou_improved_2019}. 
\begin{defn}[Smooth approximation of the indicator function] \label{def:smooth-indicator}
    Let $g : \R^d \to \R$ be $\mathcal{F}$-measurable and the set $A = \{ x \in \R^d: g(x) \leq 0\}$. Let $\sigma > 0$ be smoothing parameter, and $\mu = \mu(\sigma) > 0$ be another parameter depending on $\sigma$ such that 
\begin{align} \label{cond:sigma-mu}
\lim_{\sigma\to 0}\mu(\sigma)  = 0 \ \text{ and }  \lim_{\sigma\to 0}\frac{e^{-\mu(\sigma)/\sigma}}{\sigma} = 0.
 \end{align}
    We define the smooth approximation of $\mathds{1}_A$ as the function $F_\sigma : \R^d \to \R^+$ with
\[F_\sigma(x) = \frac{1}{1+\exp\left(\frac{-\mu + g
(x)}{\sigma}\right)} =\frac{1}{2}\left(1 + \tanh\left(-\frac{-\mu + g
(x)}{2\sigma}\right) \right),\quad x \in \R^d.\] 
\end{defn}
Note that the assumptions on $\sigma$ and $\mu$ imply that $\lim_{\sigma\to0}e^{-\mu/\sigma} = 0$. For instance, they are satisfied if $\mu=C\sqrt{\sigma}$ for some $C > 0$. Moreover, they ensure the pointwise convergence of $F_\sigma$ to $\mathds{1}_A$ as $\sigma \to 0$. Specifically, for any $x \in \R^d$ such that $g(x) \neq 0$,
\[\lim_{\sigma\to 0} F_\sigma(x)  =  \lim_{\sigma\to0} \frac{1}{1+e^{g(x)/\sigma}} = \begin{cases}
    1 &\text{ if } g(x) < 0
    \\0 &\text{ if } g(x) > 0,
\end{cases}\]
and for $g(x) = 0$, it is immediate that $\lim_{\sigma\to0} F_\sigma(x) = \lim_{\sigma\to0} (1 + e^{-\mu/\sigma})^{-1} = 1$.

 For any $\sigma > 0$ and the associated $\mu > 0$, the smoothed optimal importance sampling distribution is defined as 
\begin{align}\label{eq:smoothed-IS-dist}
\rhosig = \frac{F_\sigma \rho}{p_\sigma} \text{ with } p_\sigma = \int_{\R^d} F_\sigma(x) \rho(x)\,\d x = \E[F_\sigma(Z)],
\end{align}
and its importance sampling estimate of $p$ writes
\[\hat p_{\rho^*_\sigma} :=  \frac{1}{N}\sum_{i=1}^N\frac{ \mathds{1}_A(X_i)\rho(X_i)}{\rhosig(X_i)}, \ (X_i)_{i=1,\ldots,N} \overset{\text{i.i.d.}}{\sim} \rhosig.\]
We note that from the pointwise convergence of $F_\sigma$ to $\mathds{1}_A$ and the fact that $F_\sigma(x)\rho(x) \leq \rho(x)$ uniformly in $\sigma$, we obtain by Lebesgue's dominated convergence theorem that \\$\lim_{\sigma\to0}\lvert p_\sigma -~p \rvert =~0$. As such, we deduce a pointwise convergence of the PDF $\rho^*_\sigma$ towards $\rho^*$ as $\sigma \to 0$. This suggests that the nRMSE of $\hat p_{\rho^*_\sigma}$ should be close to that of $\hat p_{\rhoopt}$ (which itself is $0$) when the smoothing parameter $\sigma$ is small.  In Section~\ref{subsec:IS-error}, we will confirm this intuition by establishing an error bound on the nRSME of $\hat p_{\rho^*_\sigma}$.
\subsection{Different Langevin dynamics} \label{sec:Langevin-dynamics}
We aim to employ the affine-invariant interacting Langevin dynamics as an MCIS scheme for rare event estimation, by targeting the smoothed optimal distribution. We start by describing some generalities on the Langevin dynamics, before extending the discussion on the more advanced version with interaction. Then, we discuss the time-discretized version of the dynamics.
\subsubsection{Overdamped Langevin dynamics}
Let $(X_t)_{t\geq 0}$ be a time-dependent stochastic process in the state space $\R^d$, defined as solution of the following stochastic differential equation called the overdamped Langevin dynamic (often just called the Langevin dynamic):
\begin{align} \label{eq:langevin-dynamics}
\d X_t = - \nabla V(X_t) \d t + \sqrt{2} \ \d W_t, \ X_0 \sim \pi_0,
\end{align}
where $V : \R^d \to \R_+$ is called the potential function, $W_t$ is a $d$-dimensional Brownian motion and $\pi_0$ is the initial distribution.
For any $t\geq 0$, let $\pi(\cdot,t)$ be the probability density function of $X_t$ and denote by $\pi_*$ the Gibbs distribution
\[\pi_*(x) = \frac{e^{-V(x)}}{\int_{\R^d} e^{-V(z)} \, \d z},\quad x\in\R^d\,.\]
Under certain conditions the evolution of the distribution $\pi(\cdot,t)$ approaches the Gibbs distribution $\pi_*$. In the sequel we denote for $a=1,2$ and a function $\pi: \R^d \mapsto \R_+$, 
\[L^{a}(\R^d; \pi) = \left\{ f : \R^d \to \R \text{ such that } \int_{\R^d} f^a(x)\pi(x) \, \d x < +\infty \right\},\]
as well as $L^{a}_{+}(\R^d; \pi) = L^{a}(\R^d; \pi) \cap \{f:\R^d \to \R_+\}$.
We recall the following definition and result, collected in~\cite{pavliotis_14}.
\begin{defn}[{Poincaré inequality~\cite[Theorem 4.3]{pavliotis_14}}]
Let $V$ be a potential function and $\pi_*$ be its corresponding Gibbs distribution. We say that the potential $V$ or the Gibbs distribution $\pi_*$ satisfies the Poincaré inequality with constant $\lambda > 0$ if for every $f \in C^1(\R^d, \R) \cap L^2(\R^d; \pi_*)$ such that $\int_{\R^d} f(x)\pi_*(x) \, \d x = 0$, 
\[ \lambda \lVert f \rVert^2_{L^2(\R^d; \pi_*)} \leq \lVert \grad f \rVert^2_{L^2(\R^d; \pi_*)}.\]
\end{defn}
Assuming that $V$ satisfies the Poincaré inequality, one obtains the following exponential convergence result towards the stationary distribution $\pi_\ast$.
\begin{thm}[{\cite[Theorem 4.4]{pavliotis_14}}] \label{thm:gibbs-expfast}
    Assume that $\pi(\cdot,0) \in L_+^2(\R^{d}; \pi_*^{-1})$. If the potential~$V$ satisfies a Poincaré inequality with constant $\lambda > 0$, then $\pi(\cdot,t)$ converges to the Gibbs distribution exponentially fast in the following sense:
\[\lVert \pi(\cdot,t) - \pi_*\rVert_{L_+^2(\R^d;\pi_*^{-1})} \leq e^{-\lambda t}\lVert{\pi(\cdot, 0) - \pi_* \rVert}_{L_+^2(\R^d;\pi_*^{-1})}.\]    
\end{thm}
The Bakry-Emery criterion~(see e.g.~\cite[Theorem 4.4]{pavliotis_14}) provides a sufficient condition for $V$ to satisfy a Poincaré inequality. This criterion states that if $V\in C^2(\R^d)$ is strongly convex, that is, there exists some $\alpha > 0$ such that for any $ x \in \R^d, \ \Hess  (V)(x) \geq \alpha I$, then~$V$ satisfies a Poincaré inequality with constant $\alpha$. 
In the above, the term $\Hess(f)$ is the Hessian of a function $f\in C^2(\R^d,\R)$ and for any square matrices $A$ and $B$, $A \geq B$ means that $A - B$ is positive-definite. 

In other words, the overdamped Langevin dynamics allows for the sampling from the Gibbs density $\pi_* \propto e^{-V}$ as $t \to +\infty$.

\subsubsection{Affine invariant interacting Langevin dynamics~(ALDI)} \label{sec:ALDI}
The ALDI is an interacting particle based version of the Langevin dynamics introduced in~\cite{nusken2019noteinteractinglangevindiffusions,inigonuskenreich20}. The authors proved desirable theoretical properties on ALDI which we summarize below. However, the bias control on the discretized (in time) ALDI has not been studied so far, and an efficient time-stepping scheme is still subject to research. 
Despite this, ALDI will be the version we consider to propose a new algorithm because of its advantageous affine-invariance property. In this context, affine-invariance means that the dynamic of the particles are invariant w.r.t.\ an affine change of coordinates. We refer to~\cite[Definition 2]{inigonuskenreich20} for the exact definition. The advantage of the affine-invariance property for samplers has been studied both numerically \cite{CF2010, GW2010} and theoretically \cite{RS2022}. In particular, this feature is beneficial for the sampling of anisotropic target distributions, as the algorithm adapts to the anisotropy, leading to improved mixing.

Let $(X_t^{i})_{i=1,\ldots,M}$ be an ensemble of $M$ particles, each particle $X_t^{i} \in \R^d$ obeying some dynamics for $t\geq 0$, and we write $\bX_t$ for the $d \times M$ matrix where
\[\bX_t = \left( X_t^1,\ldots, X_t^M \right).\]
At each time $t$, ALDI builds on the overdamped Langevin dynamic~\eqref{eq:langevin-dynamics} by employing the empirical covariance matrix of the $M$ particles to account for the interaction. First, we define the notations for the empirical mean~$\mu$, and the empirical covariance matrix $\Sigma$ of $(X_t^{i})_{i=1,\ldots,M}$ respectively as
\[\mu(\bX_t) = \frac{1}{M}\sum_{i=1}^M X_t^{i} \quad \text{ and }\quad \Sigma(\bX_t) = \frac{1}{M}\sum_{i=1}^M (X_t^i - \mu(\bX_t))(X_t^i - \mu(\bX_t))^\top,\]
as well as the generalized square root of the empirical covariance matrix defined in~\cite{inigonuskenreich20} as
\[ \Sigma^{1/2}(\bX_t) = \frac{1}{\sqrt{M}} (\bX_t - \mu(\bX_t)\boldsymbol{1}_M^\top),\]
where $\boldsymbol{1}_M$ is a vector in $\R^M$ that contains only ones. The matrix $\Sigma^{1/2}(\bX_t)$ thus verifies $\Sigma(\bX_t) = \Sigma^{1/2}(\bX_t)(\Sigma^{1/2}(\bX_t))^\top$. Then, ALDI is defined as the following stochastic differential equation 
evolving each particle $i=1,\ldots,M$~\cite[Definition 5]{inigonuskenreich20}:
\begin{align}\label{eq:ALDI}
\d X_t^i = \underbrace{- \Sigma(\bX_t) \nabla V (X_t^i) \ \d t \vphantom{\frac{}{}}}_{\text{first drift term}} + \underbrace{\frac{d+1}{M}\left(X_t^i - \mu(\bX_t)\right) \d t}_{\text{second drift term}} + \underbrace{\sqrt{2}\ \Sigma^{1/2}(\bX_t) \ \d W_t^i \vphantom{\frac{}{}}}_{\text{diffusion term}} , 
\end{align}
where $V : \R^d \to \R_+$ is the potential function and $W_t^i$ are independent $M$-dimensional Brownian motions. The key difference between ALDI and the original Langevin dynamic is that the particles interact with each other through the introduction of the empirical covariance matrix as the preconditioner for the first drift term. The second drift term was introduced in~\cite{nusken2019noteinteractinglangevindiffusions} as a finite-size correction to obtain the correct long-time behavior $(t \to \infty)$ under the finite ensemble size regime. The third term is called the diffusion term, which involves the generalized square root $\Sigma^{1/2}(\bX_t)$, whose calculation does not require the inversion of the matrix $\Sigma(\bX_t)$.

We recall the definition of the Gibbs distribution on $\R^d$ as $\pi_* \propto \exp(-V)$. We define the following product measure on $\R^{d \times M}$,
\[\pi_*^{M}(\bX) := \prod_{i=1}^M \pi_*(X^i), \ \bX = (X^1,\ldots, X^M)\in \R^{d \times M},\]
and briefly present the ergodicity result for ALDI in~\cite{inigonuskenreich20}. 
This result requires the following assumptions on the potential $V$ which are stronger than the Bakry-Emery criterion. In the sequel, we write $\lVert \cdot \rVert$ for the euclidean norm on $\R^d$.
\begin{assumption}[{\cite[Assumption 9]{inigonuskenreich20}}] \label{assum:ALDI}
    The potential $V$ belongs to $C^2(\R^d, \R) \cap L^1(\R^d;\pi_*)$, and there exist a compact set $K \subset \R^d$ and constants $\alpha, L > 0$ such that for all $x \in \R^d \setminus K$,
    \[\alpha\lVert x \rVert^2 \leq V(x) \leq L \lVert x \rVert^2, \ \alpha\lVert x \rVert \leq \lVert \nabla V(x) \rVert \leq L\lVert x \rVert\ \text{ and } \ \alpha I \leq  \Hess(V)(x)  \leq L I. \]
\end{assumption}

Under these assumptions, the authors in~\cite{inigonuskenreich20} derived the following proposition concerning the ergodicity of ALDI.
\begin{prop}[{\cite[Proposition 10, Proposition 11]{inigonuskenreich20}}]~\label{prop:non-degen-ergo}
    Suppose Assumption~\ref{assum:ALDI} is in place and let $M > d$. 
    Then~\eqref{eq:ALDI} admits a unique global strong solution.
    If in addition $M > d +1$, then the time-marginal PDF of $\bX_t$, $\pi_t^M$ satisfies $\pi_t^M \to \pi_*^M$ as $t \to +\infty$ in total variation distance.
\end{prop}  

Note that in Proposition~\ref{prop:non-degen-ergo}, $M > d+1$ is required to ensure that ALDI has the correct invariant measure $\pi_*^M$. If $M \leq d$, then $\Sigma(\bX_t)$ is singular for all $t \geq 0$, and then $(X_t^i)$ stays in the linear subspace of $\R^d$ spanned by the $M$ initial particles $(X_0^i)$, which means that for any $t \geq 0$ and $i \leq M$, there exists $a_t^i \in \R^M$ such that 
\[ X_t^i = \sum_{j=1}^M a_t^{i}(j) X_0^j.\]
This is called the linear subspace property of ALDI. In practice, when the dimension $d$ is high, one may still wish to run ALDI with $M \leq d$ particles to save computational costs. This constitutes the motivation to consider the regularization of the covariance matrix in ALDI, first proposed in~\cite{nusken2019noteinteractinglangevindiffusions}, by replacing $\Sigma(\bX)$ by
\[\Sigma_\gamma(\bX) = (1-\gamma)\Sigma(\bX) + \gamma I,\]
for some $ 0<\gamma \leq 1$. This matrix is always positive-definite. When $\gamma = 1$, the $M$ particles evolve independently according to the overdamped Langevin dynamic~\eqref{eq:langevin-dynamics}.
The ultimate form of ALDI we will implement numerically is then the following, with an appropriate modification of the second drift term to guarantee the same convergence properties as ALDI~\cite{nusken2019noteinteractinglangevindiffusions}:
for each particle $i=1,\ldots,M$,
\begin{align}\label{eq:ALDI-regul}
\d X_t^i = - \Sigma_\gamma(\bX_t) \nabla V (X_t^i) \ \d t + (1-\gamma)\frac{d+1}{M}\left(X_t^i - \mu(\bX_t)\right) \d t + \sqrt{2}\ \Sigma_\gamma^{1/2}(\bX_t) \ \d W_t^i\,.
\end{align}
\subsubsection{Time-discretization: ULA and ALDI} \label{sec:ULA-ALDI}
For the numerical implementation of ALDI~\eqref{eq:ALDI-regul}, one needs to discretize the dynamics in time, usually through the Euler-Maruyama scheme. When $\gamma = 1$ in~\eqref{eq:ALDI-regul} we recover the overdamped Langevin dynamic, and the time-discretized version is then called the unadjusted Langevin algorithm~(ULA). In other words, let $\epsilon > 0$ be the time-step size, then ULA represents the following discrete stochastic differential equation: let $X_0 \sim \pi_0$ an initial distribution, and for $k \in \N_*$,
\begin{align}\label{eq:ULA}
 X_{k+1} = X_k - \epsilon \nabla V(X_k) + \sqrt{2\epsilon}Z_k,
\end{align}
where $Z_k \sim N(0,I)$ independent of $X_k$, and $V$ is the potential as in~\eqref{eq:langevin-dynamics}. 
Now let $\pi_k$ be the distribution of ULA at the $k$-th iteration. Recall that the Gibbs distribution is defined as $\pi_* \propto e^{-V}$. The ULA has a unique invariant measure as $k \to \infty$ if the following condition is satisfied~\cite[Lemma 1]{Wibisono2018SamplingAO}:
$V$ is $\alpha$-strongly concave and $L$-smooth for some $L \geq \alpha > 0 $, i.e.\ for any $x \in \R^d$,
\begin{equation}\label{cond:log-cc-log-smooth}
    \alpha I \leq  \Hess(V)(x) \leq L I.
\end{equation}
In addition, we have a control of bias of order $1/2$~\cite[Theorem 1]{dalalyan17a}:
\begin{thm}\label{thm:bias-ula}
    Assume~\eqref{cond:log-cc-log-smooth} and $\epsilon \leq 2/(\alpha + L)$, then for $k \geq 0$
    \[W_2(\pi_k, \pi_*) \leq (1-\alpha\epsilon)^k W_2(\pi_0,\pi_*) +    1.82(L/\alpha)(\epsilon d)^{1/2}.\]
\end{thm}
This result shows the central role the time-step size $\epsilon$ plays on the sampling bias on ULA. The condition $\epsilon \leq 2/(\alpha+ L)$ signifies that the maximum time-step size for Theorem~\ref{thm:bias-ula} to hold decreases with the concavity and the smoothness of the potential. This result relies on an assumption that is weaker than Assumption~\ref{assum:ALDI} used for ALDI. 

The time-discretized version of ALDI~\eqref{eq:ALDI-regul} with time-step size $\epsilon > 0$ writes for $k \in \N_*$, 
and for each particle $i=1,\ldots,M$,
\begin{multline}\label{eq:discrete-ALDI-fixedsize}
    X_{{k+1}}^i = \left[ - \Sigma_\gamma(\bX_{k}) \nabla V (X_{k}^i) + (1-\gamma)\frac{d+1}{M}\left(X_{k}^i - \mu(\bX_{k})\right)\right] \epsilon
    + \sqrt{2 \epsilon}\ \Sigma_\gamma^{1/2}(\bX_{k})Z_{k}^i
\end{multline}
where $(Z_{k^i})_{i=1,\ldots,M}$ are i.i.d.\ $d$-dimensional standard Gaussian vectors and independent of $(X_{k}^i)_{i=1,\ldots,M}$. The time-discretized version of ALDI is still called ALDI. A result along the lines of Theorem~\ref{thm:bias-ula} for ALDI yet remains an open question. A practical workaround consists in introducing a Metropolization step as in~\cite{sprungk2025} to circumvent the sampling bias due to time discretization. 

\subsection{Contribution and organization of the work}
In Section~\ref{sec:aldi-IS}, we propose a new MCIS scheme for rare event estimation based on ALDI~\eqref{eq:ALDI-regul} introduced in~\cite{nusken2019noteinteractinglangevindiffusions,inigonuskenreich20}, termed affine invariant Langevin dynamics for importance sampling~(ALDI-IS). We present several practical techniques proposed for implementation, including sequential rare event sampling in this context and the coupling of ALDI with a clustering algorithm.

In Section~\ref{sec:MCIS-theoretical}, we first establish a non-asymptotic error bound on the nRMSE of the importance sampling estimator based on $\rhosig$ (Theorem~\ref{thm:error-bound}). While a smaller smoothing parameter $\sigma$ reduces the importance sampling error, it simultaneously makes sampling from $\rhosig$ increasingly difficult. We will analyze this trade-off by quantifying the sampling bias of the ULA, a special case of ALDI (Proposition~\ref{prop:all-conditions}). More precisely, we observe a trade-off in the choice of $\sigma$ between the ULA step and the IS step, since increasing $\sigma$ reduces the bias on the former step. The proofs of this section are collected in~\ref{sec:proofs}.

In Section~\ref{sec:num}, we compare the performance of ALDI-IS with recent state-of-the-art estimators on standard rare event estimation benchmarks.

\section{A new gradient-based MCIS scheme: ALDI-IS} 
\label{sec:aldi-IS}
We propose a new gradient-based MCIS scheme using the ALDI, in line with the works of~\cite{wagner22,PAPAKONSTANTINOU2023,althaus2024}, and discuss our numerical strategies to be adopted in the context of rare event estimation. As outlined in Section~\ref{subs:Intro-IS}, MCIS is structured into three distinct steps. We now describe the specific implementation of ALDI within this framework, addressing each step in turn.

\subsection{Step 1: 
ALDI for rare event estimation}
First, we show how the ALDI can be employed as an MCIS algorithm for rare event estimation. We have seen in Proposition~\ref{prop:non-degen-ergo} (and Theorem~\ref{thm:gibbs-expfast} for $\gamma = 1$) that if a potential $V$ satisfies Assumption~\ref{assum:ALDI}, then it is possible to sample from the corresponding Gibbs distribution. 
In our context of rare event estimation, to sample from $\rhosig$, we define the following potential.
\begin{defn}\label{def:rare-event-potential}
    The rare event potential $V_\sigma:\R^d\to\R$ is defined by
    \[V_\sigma(x) = - \log F_\sigma(x) + \frac{1}{2}\lVert x \rVert^2,\quad x\in\R^d\,.\]   
\end{defn}
The corresponding Gibbs density for $V_\sigma$ coincides with the smoothed optimal importance sampling density, since for all $x \in \R^d$,
\[e^{-V_{\sigma}(x)} = F_\sigma(x) e^{-\lVert x \rVert^2/2} = F_\sigma(x)\rho(x) \propto \rhosig(x). \]
Hence, 
we can apply the ALDI with $V_\sigma$ to sample from the smoothed optimal density~$\rhosig$ provided that $V_\sigma$ satisfies Assumption~\ref{assum:ALDI}. We give a quick intuition as to why the assumption can be satisfied.
Since $\lVert x \rVert^2$ is a quadratic potential, $V_\sigma$ can be viewed as a perturbation of a quadratic potential by $-\log F_\sigma$, so that we reasonably expect that this assumption can be satisfied provided that $ - \log F_\sigma$ (thus $g$) does not grow too fast as $\lVert x \rVert \to +\infty$. We will derive a more precise quantification in Section~\ref{sec:MCIS-theoretical}.

For the next subsections, we describe additional numerical strategies to adapt ALDI for rare event estimation. These are mainly guided by practical experience, and we emphasize that there is no theoretical verification of our proposed heuristics.
\subsubsection{Adaptive time-step size } \label{sec:adaptive-timestep}
As discussed, in practice ALDI~\eqref{eq:ALDI-regul} needs to be discretized in time. We will show in Section~\ref{subs:ULA-numerics} that small smoothing parameters $\sigma$ for $\rho_\sigma$ have a detrimental effect on the sampling bias on ULA, and one needs a smaller time-step size in order for ULA to be stable. To counteract the effect of small smoothing parameters, we choose the Euler-Maruyama scheme with an adaptive time-step size similar to~\cite[Section 4.4]{reich21}, in which the step sizes are inversely proportional to the highest norm of the drift term among the particles. The discrete time $(t_k)_{k\in \N}$ is defined as $t_0 = 0$, and for all $k \in \N$, $ t_{k+1} = t_{k} + \Delta t_k$, where
\begin{align} \label{eq:adaptive-time-step}
 \Delta t_k := 0.1\cdot \left(\max_{i=1,\ldots, M}\left\{\left\lVert- \Sigma_\gamma(\bX_{t_k}) \nabla V (X_{t_k}^i)  + (1-\gamma)\frac{d+1}{M}\left(X^i_{t_k} - \mu(\bX_{t_k})\right)\right\rVert \right\}\right)^{-1}\,.
\end{align}
The discretized scheme of the ALDI~\eqref{eq:ALDI-regul} is then given by:
let $(X_0^i)_{i=1,\ldots M}$ i.i.d.\ according to~$\rho$, and for $k \in \N$, for each particle $i=1,\ldots,M$,
\begin{multline}\label{eq:discrete-ALDI}
    X_{t_{k+1}}^i = \left[ - \Sigma_\gamma(\bX_{t_k}) \nabla V (X_{t_k}^i) + (1-\gamma)\frac{d+1}{M}\left(X_{t_k}^i - \mu(\bX_{t_k})\right)\right] \Delta t_k\\
    + \sqrt{2 \Delta t_k}\ \Sigma_\gamma^{1/2}(\bX_{t_k})Z_{t_k}^i
\end{multline}
where $(Z_{t_k^i})_{i=1,\ldots,M}$ are i.i.d.\ $d$-dimensional standard Gaussian vectors and independent of $(X_{t_k}^i)_{i=1,\ldots,M}$. 

\subsubsection{Adaptive stopping criterion} \label{sec:stopping-criterion}
We have seen the ergodic property for ALDI that states that as $t \to \infty$, the distribution of the particles tends to the correct invariant distribution $\pi_*$ which we chose as $\rhosig$. In practice, the dynamic has to be stopped at a finite time. We propose the following stopping criterion when using the rare event potential $V_\sigma$, based on the cumulative sum up to the current iteration of the norm of $\nabla V_\sigma$. We define $V^{\cumu}_{-1} = 0$, and for any iteration $k \in \N$, 
\[ V^{\cumu}_k = \frac{1}{k+1}\frac{1}{M}\sum_{i=1}^M \left( \lVert \nabla V_\sigma(X_{t_k}^i)\rVert^2 + \lVert X_{t_k}^i  \rVert^2\right) + \frac{k}{k+1}V^{\cumu}_{k-1}.\]
Choose some $ 0< \varepsilon^{\cumu} < 1$ and $k_{\min} > 0$. The dynamics is stopped after $t_k$ where $k \geq k_{\min}$ is the first iteration such that
\begin{align}\label{eq:stopping-criteria}
    \frac{\lvert V^{\cumu}_k - V^{\cumu}_{k-1}\rvert}{V_k^{\cumu}} \leq \varepsilon^{\cumu}.
\end{align}
We will denote the stochastic stopping time as $T$. The stopping criterion is based on the heuristic that the integration over the norm of the two drift terms, $\lVert \nabla V_\sigma(\cdot)\rVert^2 + \lVert \cdot  \rVert^2$ w.r.t.\ the invariant distribution, should be approximated by the ergodic theorem and thus stop evolving in time.
Note that we require the introduction of the minimum number of iterations $k_{\min}$, since the evolution of the gradient might be small in the first few iterations.

\subsubsection{Sequential rare event sampling}
Recall that the rare event set is defined as $A = \{ x \in \R^d : g(x) \leq 0\}$. In numerical experiments, $A$ can be approached sequentially, which is a classical approach used in most rare event simulation algorithms. Two of the most popular approaches are:
\begin{enumerate}
    \item Introduce a decreasing sequence of $Q$ positive numbers $(q_j)_{j=1,\ldots Q}$ and consider $A_j = \{ x \in \R^d : g(x) - q_j \leq 0\}$. Then run the sampling scheme for each $A_j$ by using samples from the previous iteration $j-1$;
    \item Gradually decrease the smoothness of $F_\sigma$ by decreasing $\sigma$ iteratively.
\end{enumerate}
The first approach is used for instance in subset simulation approaches~\cite{au_estimation_2001,cerou_sequential_2012,elfverson2023adaptivemultilevelsubsetsimulation} and the adaptive importance sampling scheme cross-entropy~\cite{rubinstein_cross_2004}. Fixed-level subset simulation considers a fixed sequence $(q_j)$ beforehand, whereas cross-entropy and adaptive subset simulation decrease $q_j$ adaptively. For our proposed scheme, we will fix $(q_j)$, and replace $g$ by $g - q_j$ in the definition of $F_\sigma$. Thus, we modify the potential $V_\sigma$ to be used in ALDI.

The second approach is used in the improved cross-entropy algorithms~\cite{papaioannou_improved_2019,ELMASRI2021,uribe2021,wagner22} and the adaptive decrease in $\sigma$ is often calculated based on the coefficient of variation of the likelihood weights of two subsequent PDF's. However, for Langevin dynamics, the PDF of the samples has to be estimated unlike improved cross-entropy algorithms, and the smoothness parameter $\sigma$ affects the bias via Corollary~\ref{cor:bias-ula-sigma}. We instead decide to fix a $\sigma$ as in~\cite{ehre2024steinvariationalrareevent} and consider the first approach of sequentially approaching the rare event set $A$. 

Moreover, we also consider a decreasing sequence of regularization parameters $(\gamma_j)$. The underlying intuition is that when $q_j$ is high, then $A_j$ is a larger subset of $\R^d$. Specifically, ULA without interaction should work well and allow for exploration of the whole $\R^d$. As $q_j$ decreases, the particles will concentrate towards the subset $A_j$, which is closer to $A$ (in the inclusion sense). If one makes the educated guess that $A$ is embedded in a linear subspace of $\R^d$ with lower dimension in the rare event estimation, then $A_j$ should be as well for small enough $q_j$, so one can use $\Sigma(\bX_t)$ even when $M < d$, since the subspace property can be exploited intelligently if the previous evolution of the particles has led towards the correct linear subspace using the previous larger subsets. A large $\gamma$ allows for the exploration of the whole space $\R^d$ to reach the correct subsets, whereas a small $\gamma$ exploits the subspace spanned by the particles to accelerate the convergence once the correct subsets are found. The tuning of $(\gamma_j)$ needs to take into account the trade-off between exploration and exploitation. In short, tuning a higher $\gamma$ implies more exploration which is suitable in the beginning of the algorithm where the particles should ``search for'' the rough position of rare event set, meanwhile tuning a lower $\gamma$ implies more exploitation which is suitable for later stages of the algorithm once the particles are close to the rare event set.

\subsubsection{Gradient of mean particles with clustering algorithm}
In practical applications, the limit state function $g$ is often unimodal or multimodal w.r.t.~the minimum. In such cases, the set $A$ is termed a unimodal or a multimodal hyperspace in $\R^d$. Consequently, we would like to design low interactions between particles belonging to different modes, while strong interactions between particles within the same mode. To deal with multi-modality, our strategy uses mixture models for the importance sampling auxiliary distribution such as Gaussian mixtures and von-Mises Fisher Nakagami mixture models~\cite{papaioannou_improved_2019,geyer_cross_2019}, and employs ALDI in a similar fashion as using gradient-based optimization to search for failure points near the modes~\cite{chiron_failure_2023}.

Since the gradient of each particle is evaluated in ALDI, the particles naturally evolve towards the different modes in a multimodal hyperspace if the initial spread of the particles is large enough. 
However, around the same mode, the gradient of the particles tends to point towards the direction of the maximum a posteriori in the mode, so if one successfully identifies that particles within a cluster are moving towards the same mode, one can potentially save gradient calls by only evaluating a representative gradient instead of that of each particle in the cluster.
To this effect, we propose to couple a clustering algorithm in ALDI to identify clusters of particles, and for each particle in the same cluster, only the gradient of the mean particle is to be evaluated. 

In our practical implementation, we use the DBSCAN algorithm~\cite{dbscan1996} implemented in the Statistics and Machine Learning Toolbox of Matlab. Our clustering algorithm starts after $10$ burn-in iterations to allow the particles to naturally form clusters beforehand. For the parameters of DBSCAN in the Matlab implementation, we use the heuristic that the epsilon neighbourhood should be $d/2$. The minimum numbers of neighbours required for core points are $5$ when the total number of particles are at least $50$. We note that if ALDI is run with a number of particles $M$ lower than $50$, one should decrease the number of neighbours required accordingly. Specifically for ALDI, it is heuristically better to apply overly strict parameters for clustering compared to overly convenient, since if the parameters are too strict, ALDI-DBSCAN coincides with ALDI with full gradient calls and provides accurate results albeit with additional gradient calls. However, if the parameters are too lenient, one runs the risk of the particles which are supposed to evolve towards two separate modes being identified as one single cluster, for instance. Since only the gradient of the mean particle is calculated, the particles would not evolve towards the correct modes: each particle will evolve towards the mode of the mean particle movement. For a better illustration, we refer to Figure~\ref{fig:Fourbranches} in Section~\ref{sec:four-branches}. In addition, the frequency of application of DBSCAN on the particles is usually once every 10 iterations, and we require 10 burn-in iterations. 
\subsection{Step 2: von-Mises Fisher Nakagami mixture model}~\label{subsection:vMFNM} 
After the final iteration of ALDI is achieved at time $T$, we collect the $M$ samples $(X_T^i)_{i=1,\ldots,M}$ which are approximately distributed according to 
the smoothed optimal density $\rhosig$. In order to use $(X_T^i)_{i=1,\ldots,M}$ in rare event estimation, we first fit a parametric distribution model onto the particles, then generate new samples from the fitted parametric distribution to estimate $p$ by importance sampling to guarantee the unbiasedness. This approach introduces additional parametric fitting error which will be explored in a future study. Since we desire to run ALDI with $M \leq d$ samples, the $K$-Gaussian mixtures distribution model which contains $Kd(d+3)/2 + (K-1)$ scalar parameters to estimate ($d$ and $d(d+1)/2$ scalars for the mean and covariance matrix of each of the $K$ mixtures, and $K-1$ mixture coefficients) can be too high to be tuned properly. We instead focus on the von-Mises Fisher Nakagami mixture (vMFNM) model, introduced for rare event estimation in~\cite{wang_cross-entropy-based_2016,papaioannou_improved_2019} to approximate a high-dimensional Gaussian mixture with less parameters. With $K=1$ mixture, the parameters of the von-Mises Fisher Nakagami (vMFN) distribution are given by
\[\nu = (\mu,\kappa,m,S), \text{ where } \mu \in \mathcal{S}^{d-1} := \{x \in \R^d: \lVert x \rVert = 1 \}, \ \kappa \geq 0, \ m \geq 0.5 \text{ and } S > 0.\]
Using the polar coordinates on $\R^d = \R_{+} \times \mathcal{S}^{d-1}$, where we write for any $x \in \R^d \setminus \{0_{\R^d}\}$,
\[ r = \lVert x \rVert > 0 \text{ and } \omega = \frac{x}{\lVert x \rVert} \in \mathcal{S}^{d-1},\]
so that $x = r \times \omega$ (and $0_{\R^d} = 0_{\R_+} \times \omega_0$ for any $w_0 \in \mathcal{S}^{d-1}$),  the PDF of the vMFN distribution with parameters $\nu$, $\rho_\nu^{\vMFN}$, can be decomposed as
\[\rho_\nu^{\vMFN}(r,\omega) = \rho_{m,S}^{\text{N}}(r)\rho_{\mu,\kappa}^{\vMF}(\omega),\]
where $\rho_{m,S}^{\text{N}}$ is the PDF of the Nakagami distribution with parameters $(m, S)$~\cite{NAKAGAMI19603}. This distribution is defined as 
\[ \rho_{m,S}^{\text{N}}(r) = \frac{2m^m}{\Gamma(m) S^m} r^{2m-1} \exp\left(-\frac{m}{S}r^2\right), \ r \in \R_+, \]
where $\Gamma$ is the Gamma function which models the radius $\lVert x \rVert$ of $x \in \R^d$. The parameter $m$ is called the shape parameter whereas the parameter $S$ is called the spread parameter. The distribution $\rho_{\mu,\kappa}^{\vMF}$ is the von-Mises Fisher distribution~\cite{wang_cross-entropy-based_2016} with parameters $(\mu,\kappa)$ defined by
\[\rho_{\mu,\kappa}^{\vMF}(\omega) = \frac{\kappa^{d/2 - 1}}{(2\pi)^{d/2}I_{d/2 - 1}(\kappa)} e^{\kappa \mu^\top\omega}, \ \omega \in \mathcal{S}^{d-1},\]
with $I_\alpha$ the modified Bessel function of the first kind of order $\alpha$. The parameter $\mu$ is called the mean direction parameter, and the parameter $\kappa$ the concentration parameter. This distribution characterizes the distribution of the probability mass on $\mathcal{S}^{d-1}$, with concentration around the direction of $\mu$.
With $K$ mixtures, the vMFNM model with parameters $\boldsymbol{\nu} = (\nu_1,\ldots,\nu_K)$ writes
\[\rho_{\boldsymbol{\nu}}^{\vMFNM}(r,\omega) = \sum_{i=1}^K \alpha_i \rho^{\vMFN}_{\nu_i}(r,\omega),\]
where $(\alpha_i)$ are $K$ weights that sum to $1$. There are $K(d+3) + (d-1)$ parameters, which is more amenable to be tuned via the expectation-maximization algorithm than the Gaussian mixture model, when one has $M \leq d$ particles from ALDI. For the details of the expectation-maximization algorithm for vMFNM model and the choice of $K$, we refer to~\cite[Section 5.1]{papaioannou_improved_2019}.
\begin{algorithm}[!htbp]
\caption{Affine invariant Langevin dynamics for importance sampling (ALDI-IS)}\label{alg:ALDI-IS}
\begin{algorithmic}
\Require \textbf{Hyperparameters, for which we have good heuristics} \\
\begin{tabular}{ll}
Smoothing parameters &$\sigma > 0, \ \mu > 0$ verifying~\eqref{cond:sigma-mu} \\
Number of rare event levels     &$Q \in \N_*$ \\
Rare event levels &$(q_j)_{j=1,\ldots Q}, q_j \geq 0 \ \forall j$ \\
     Regularization parameters &$(\gamma_j)_{j=1,\ldots, Q}, \ 0 < \gamma_j < 1 \ \forall j$ \\
     Stopping criteria for ALDI &$(\varepsilon^{\cumu}_j)_{j=1,\ldots, Q}, \ 0 < \varepsilon^{\cumu}_j \leq 1$ \\
    Number of particles for ALDI &$M$ \\
    Number of samples for IS &$N$ 
    \\
    \textbf{If DBSCAN:}\\
    Epsilon neighboorhood &$\varepsilon \approx d/2$\\
    Minimum number of neighbours &$\delta \, \approx 5$
\end{tabular} \\ \ \\
\State 1. Start with $M$ particles $(X_0^i)_{i=1,\ldots M}$ i.i.d.\ according to $\rho = N(0,I)$. \\
\State 2. Run ALDI:
\For{ $j=1,\ldots,Q$}
    \State (a) Define the smoothed approximation of the indicator function of the set $A_j = \{ x \in \R^d: g(x) - q_j \leq 0 \}$:
    \[F_\sigma(x) = \left(1+\exp\left(\frac{-\mu + g(x) - q_j
}{\sigma}\right)\right)^{-1}, \  x \in \R^d.\] 
    \State (b) Run the discrete ALDI~\eqref{eq:discrete-ALDI} on $(X_0^i)_{i=1,\ldots,M}$ with regularization parameter $\gamma = \gamma_j$ and with the potential $V_\sigma(x) = - \log F_\sigma(x) + \lVert x\rVert^2/2, \ x \in \R^d$,
    using discrete times $(t_k)$ defined in~\eqref{eq:adaptive-time-step}, until the stochastic stopping criterion~\eqref{eq:stopping-criteria} is met with $\varepsilon^{\cumu} = \varepsilon^{\cumu}_j$. This gives $M$ particles $(X_T^i)_{i=1,\ldots,M}$. 
    \\\textbf{If DBSCAN:}\\
    \quad\quad$\triangleright$ After 10 burn-in iterations, and for every tenth iteration, DBSCAN with parameters $(\varepsilon,\delta)$ is applied on $(X^i_t)_{1,\ldots,M}$ to obtain $D$ clusters $\{C_1, \ldots, C_D\}$. 
    \\\quad\quad$\triangleright$ For each cluster $k$, calculate the mean particle $X^{C_k}_t$ and its gradient of the potential $\nabla V_\sigma(X^{C_k}_t)$. 
    \\\quad\quad$\triangleright$ Assign for $i= 1,\ldots, M$,
    $\nabla V_\sigma(X_t^i) \leftarrow \nabla V_\sigma(X^{C_k}_t)$ if there exists $k= 1,\ldots, D$ such that $X_t^i \in C_k$. Otherwise (outlier particles), calculate $\nabla V_\sigma(X_t^i)$.
    \\
    \State (c) Assign $(X_0^i)_{i=1,\ldots,M} \leftarrow (X_T^i)_{i=1,\ldots,M}$.
    \EndFor
    \\
\State 3. Fit a vMFNM model on the particles $(X_T^i)_{i=1,\ldots,M}$ using the EM algorithm defined in~\cite[Section 5.1]{papaioannou_improved_2019}. This gives the distribution $\rho_{\hat{\boldsymbol{\nu}}}^{\vMFNM}$.
\\
\State 4. Importance sampling step:
\\Generate $N$ samples $(X_i)_{i=1,\ldots,N}$ i.i.d.\ according to $\rho_{\hat{\boldsymbol{\nu}}}^{\vMFNM}$ and construct the importance sampling estimator of $p$ as 
\[\hat p_{\rho_{\hat{\boldsymbol{\nu}}}^{\vMFNM}} = \frac{1}{N}\sum_{i=1}^N\frac{ \mathds{1}_A(X_i)\rho(X_i)}{\rho_{\hat{\boldsymbol{\nu}}}^{\vMFNM}(X_i)}.\]
\end{algorithmic}
\end{algorithm}
\subsection{Step 3: importance sampling estimate of rare event probability} 
Once the vMFNM model is fitted on the $M$ particles $(X_T^i)$ of ALDI, which results in the estimated parameter $\hat{\boldsymbol{\nu}}$, we obtain the distribution $\rho_{\hat{\boldsymbol{\nu}}}^{\vMFNM}$. Finally, we generate $N$ new i.i.d.~samples $(X_i)_{i=1,\ldots,N}$ according to $\rho_{\hat{\boldsymbol{\nu}}}^{\vMFNM}$ and construct the importance sampling estimator of $p$ as in~\eqref{def:IS-estimator}, i.e. 
\[\hat p_{\rho_{\hat{\boldsymbol{\nu}}}^{\vMFNM}} = \frac{1}{N}\sum_{i=1}^N\frac{ \mathds{1}_A(X_i)\rho(X_i)}{\rho_{\hat{\boldsymbol{\nu}}}^{\vMFNM}(X_i)}.\]
The complete algorithm that we name as affine invariant Langevin dynamics for importance sampling (ALDI-IS), is summarized in Algorithm~\ref{alg:ALDI-IS}, with the option of employing DBSCAN.

\section{Trade-off between importance sampling error and discretization bias with respect to the smoothing parameter}\label{sec:MCIS-theoretical}
In this section, we carry out a theoretical study on the influence of the smoothing parameter $\sigma$ occuring in the smoothed optimal density $\rhosig$. First, we study the normalized root mean square error of the importance sampling estimator if one were able to generate i.i.d.\ samples directly from $\rhosig$ and know exactly its PDF to evaluate the IS weights. Then, we check how the smoothing parameter influences the sampling bias of ULA, which can be viewed as a specific time discretization for ALDI with $\gamma = 1$. The analysis for $\gamma < 1$ is more involved and is left for future work. Similarly, the error committed while fitting the PDF of the final particles with a parametric distribution model is left to an upcoming work as it is of independent interest. 

We will observe that the error bound on the importance sampling error decreases as $\sigma$ decreases, which is an intuitive result since $\rho_\sigma^*$ tends pointwise to the optimal IS distribution $\rhoopt$ as $\sigma \to 0$. However, since $\rhoopt$ is not differentiable, it is impossible to sample from it using ULA. We will see that in fact an increasing $\sigma$ decreases the sampling bias using ULA. There is then a trade-off on the parameter $\sigma$ between the low importance sampling error and the feasibility of correctly sampling using ULA.

All the proofs in this section are collected in Section~\ref{sec:proofs}.

\subsection{Influence of smoothing on the importance sampling error} \label{subsec:IS-error}
In this section, we study the impact of the smoothing parameter $\sigma$ on the nRMSE of the IS estimator of $p$ constructed using the smoothed optimal density, i.e., the nRMSE of the following estimator:
\[\hat p_{\rho^*_\sigma} :=  \frac{1}{N}\sum_{i=1}^N\frac{ \mathds{1}_A(X_i)\rho(X_i)}{\rhosig(X_i)}, \ (X_i)_{i=1,\ldots,N} \overset{\text{i.i.d.}}{\sim} \rhosig,\]
for some $\sigma > 0$. One expects that a small $\sigma$ is beneficial since the smoothed optimal density $\rhosig$ is closer to the optimal density $\rhoopt$. We present 
our first main theoretical result. \newpage


\begin{thm}\label{thm:error-bound}
    For $N \geq 1$ and $\sigma > 0$, let $(X_i)_{i=1,\ldots,N}$ be i.i.d.\ according to $\rhosig$ and $Z \sim \rho$, as defined above. Assume that $g$ verifies $\E[g(X_1)^2] < +\infty$. Then, for any random variable $Y \sim \rho^*$,
    \begin{align*}
     &\frac{1}{p}\E\left[\left\lvert \hat p_{\rhosig} - p \right\lvert^2\right]^{1/2}
     \\&\leq \frac{6}{\sqrt{N} p} \left( e^{-\mu/\sigma}+ \P(0 < g(Z) <  2\mu) + \frac{e^{-2\mu/\sigma}}{\sigma^2} \E\left[\left\lvert g(X_1)- g(Y)\right\rvert^2 \indicatorA{X_1}\right]\right)^{1/2}
\end{align*}
\end{thm}
Two remarks are in place to take regarding the achieved error bound.
\begin{rk}\label{rk:IS-error-bound}
We discuss two terms that arise in the bound in Theorem~\ref{thm:error-bound}.
\begin{enumerate}
    \item Provided the existence of $C > 0$ such that $\sup_{0 <\sigma < C}\E\left[\left\lvert g(X_1)- g(Y)\right\rvert^2 \indicatorA{X_1}\right]$ is finite, the term $\frac{e^{-2\mu/\sigma}}{\sigma^2} \E\left[\left\lvert g(X_1)- g(Y)\right\rvert^2 \indicatorA{X_1}\right]$ decays to $0$ as $\sigma \to 0$ under the assumption that $\sigma$ and $\mu$ verify~\eqref{cond:sigma-mu}. This is for instance true if there exists $\sigma > 0$ such that $\E[g(X_1)^2] < +\infty$, due to the following lemma.
    \begin{lemma}\label{lem:sigma-implies-all}
        If~\eqref{cond:sigma-mu} holds and that there exists $\sigma > 0$ such that $\E[g(X_1)^2] < +\infty$, then 
        \begin{itemize}
            \item $\E[g(Y)^2] < +\infty$;
            \item  $\E[g(X_1)^2] < +\infty$ for any $\sigma > 0$; and
            \item $\lim_{\sigma \to 0} \E[g(X_1)^2] = \E[g(Y)^2]$.
        \end{itemize}
    \end{lemma}
    
    On the contrary, if for any $\sigma > 0$, $\E[g(X_1)^2]$ is infinite, then we can exploit some flexibility in the choice of $g$: In fact, we are ultimately interested in estimating the probability $\P(Z \in A)$ where $Z \sim \rho$ and $A = \{g \leq 0\}$, and the probability remains the same if $g$ is changed without changing the set $A$. For example, we also have $A = \{g^3 \leq 0\}$. More generally, one can compose $g$ with any transformation $h : \R \to \R$ that preserves the set $A$, i.e., such that $\{g \leq 0\} = \{h \circ g \leq 0\}$. Another example of $h$ that makes $h \circ g \in L^2(\R^d;\rhosig)$ while verifying $\P(g(Z) \leq 0) = \P(h \circ g(Z) \leq 0)$ is the truncation function $h_r$ for $r > 0$, defined as
\begin{align}\label{eq:truncation-function}
    h_r(x) = \textnormal{sign}(x)\min\{\lvert x \rvert , r\}, \ x \in \R.
\end{align}
\item Suppose that \eqref{cond:sigma-mu} is in place, then since $\mu \to 0$, the probability $\P(0 < g(Z) < 2\mu)$ tends to zero as $\sigma \to 0$. More precisely, if $g(Z)$ possesses a PDF $\rho_{g}$ which is continuous in $0$ and satisfies $\rho_{g}(0) > 0$, then by Taylor's theorem it holds
\[\P(0 < g(Z) <  2\mu) \underset{\sigma \to 0}{\sim} 2 \rho_g(0) \mu\,,\] 
implying that $\P\left(0 < g(Z) < 2\mu\right) \to 0$ as $\sigma \to 0$ at the same rate as $\mu(\sigma) \to 0$. If instead $\rho_g(0) = 0$, then we have $\P(0 < g(Z) <  2\mu) = o(\mu)$, and in some cases one can, for instance, refer to the exponential Tauberian theorem~\cite[Theorem 4.12.10 (iii)]{Bingham_Goldie_Teugels_1987} for a finer control of the rate of convergence of $\P\left(0 < g(Z) < 2\mu\right)$ to $0$.
\end{enumerate}
\end{rk}

If $g$ is $K$-Lipschitz continuous w.r.t.\ the euclidean norm, i.e., $\lvert g(x) - g(y) \rvert \leq K \lVert x - y \rVert$ for $x,y \in \R^d$, then the term $\E [ \lvert g(X_1)- g(Y)\rvert^2 \indicatorA{X_1} ]$ appearing in Theorem~\ref{thm:error-bound} can be simplified. Indeed, in this case we have
\[ \E\left[\left\lvert g(X_1)- g(Y)\right\rvert^2 \indicatorA{X_1}\right] \leq K^2 \E\left[\left\lVert X_1- Y\right\rVert^2\right]. \]
In this bound, the pair $(X_1, Y)$ has an arbitrary distribution with given marginals $X_1 \sim \rho^*_\rho$ and $Y \sim \rho^*$, and thus minimizing $\E\left[\left\lVert X_1- Y\right\rVert^2\right]$ over the joint distribution of $X_1$ and $Y$ we obtain, for this specific coupling $(X_1, Y)$,
\begin{align}\label{eq:minimize-wasserstein}
\E\left[\left\lvert g(X_1)- g(Y)\right\rvert^2 \indicatorA{X_1}\right] \leq K^2 \left(W_2(\rho^*_\sigma,\rho^*)+ \varepsilon \right) \text{ for any } \varepsilon > 0, 
\end{align}
where for two measures $\rho_1$ and $\rho_2$ on $\R^d$, $W_2(\rho_1,\rho_2)$ is the 2-Wasserstein distance, defined by
\[W_2(\rho_1,\rho_2) = \left(\inf_{\xi \in \mathcal{H}(\rho_1,\rho_2)}\int_{\R^d \times \R^d} \lVert x - y \rVert^2 \, \xi(\d x, \d y)\right)^{1/2},\]
where $\mathcal{H}(\rho_1,\rho_2)$ is the set of couplings on $\R^d \times \R^d$ with marginals $\rho_1$ and $\rho_2$. Since Theorem~\ref{thm:error-bound} allows for any choice of $Y\sim \rho^*$, and by letting $\varepsilon \to 0$ in the RHS of~\eqref{eq:minimize-wasserstein}, we get the following corollary.

\begin{corollary}\label{cor:error-bound-wass}
    In the setting of Theorem~\ref{thm:error-bound}, if in addition $g$ is $K$-Lipschitz continuous w.r.t.\ the euclidean norm on $\R^d$, then 
    \[
     \frac{1}{p}\E\left[\left\lvert \hat p_{\rhosig} - p \right\lvert^2\right]^{1/2} 
          \leq \frac{6}{\sqrt{N} p} \left( e^{-\mu/\sigma}+ \P(0 < g(Z) <  2\mu) + \frac{K^2 e^{-2\mu/\sigma}}{\sigma^2}  W^2_2(\rhosig, \rhoopt)\right)^{1/2}.
    \]
\end{corollary}
\subsection{Influence of smoothing on the discretization bias of ULA}~\label{subsection:ULA}
We now turn our focus to the ease of sampling from the smoothed optimal distribution. For simplicity, we work with the theoretically well-established ULA. 
More precisely, we quantify the effect of the smoothing parameter $\sigma$ of the smoothed distribution $\rhosig$ on the bias of ULA w.r.t.\ the time discretization step size, and show that a larger smoothing parameter in fact decreases the bias, making a smoother optimal distribution (higher $\sigma$) easier to sample using ULA. This leads to a compromise to be made between the following results and Theorem~\ref{thm:error-bound} in which the importance sampling error bound decreases as $\sigma$ decreases. This effect mimics the natural intuition, since as $\sigma \to 0$, $\rhosig \to \rhoopt$ which is not a smooth distribution, so decreasing $\sigma$ decreases the `smoothness' of $\rhosig$ needed for gradient-based sampling schemes.

\subsubsection{Sampling from smoothed optimal density using ULA}
We have discussed the theoretical properties of ULA in Section~\ref{sec:ULA-ALDI}. The aim of this section is to verify that with some additional assumptions on $g$, the smoothed optimal IS distribution $\rhosig$ can be formulated as a Gibbs distribution which satisfies condition~\eqref{cond:log-cc-log-smooth}, so that particles can be sampled from $\rhosig$ with ULA and its bias controlled via Theorem~\ref{thm:bias-ula}. 
Recall that the rare event potential~\eqref{def:rare-event-potential} defined for all $x \in \R^d$ as
\[V_\sigma(x) = - \log F_\sigma(x) + \frac{1}{2}\lVert x \rVert^2,\]
that allows the sampling of particles distributed according to the smoothed optimal density $\rhosig$. We aim to quantify the impact of the smoothness parameter on the bias of ULA using Theorem~\ref{thm:bias-ula}.

The main result of this subsection is the following proposition, where without loss of generality we announce for bounded $g$ (see~\eqref{eq:truncation-function} of Remark~\ref{rk:IS-error-bound}).
\begin{prop} \label{prop:all-conditions}
    Let $g$ of class $C^2(\R^d, \R)$ and bounded by some $r > 0$. Assume that there exist $K > 0$ and $G > 0$ such that
    \begin{enumerate}
        \item the Hessian of $g$ verifies $\inf_{x \in \R^d} \lambda_{\min}(\Hess(g)(x)) > - \infty$, where for any $x\in\R^d$, $\lambda_{\min}(M)$ is the smallest eigenvalue of a square matrix $M$,
        \item $g$ is $K$-Lipschitz continuous w.r.t.\ the euclidean norm on $\R^d$, i.e.\ for any  $x , y \in \R^d$,
        \[\lvert g(x)-g(y)\rvert \leq K\lVert x - y \rVert,\]
        \item $g$ is $G$-smooth ($\grad g$ is $G$-Lipschitz continuous) w.r.t.\ the euclidean norm on $\R^d$, i.e.\ for any $x,y \in \R^d$,
        \[\lVert \nabla g(x) - \nabla g(y) \rVert \leq G \lVert x - y \rVert.\]
    \end{enumerate}
Define for any $\sigma > 0$,
\[H_\sigma^g := \inf_{x \in \R^d}(1-F_\sigma(x))\lambda_{\min}(\Hess(g)(x)).\]
Then we have $H_\sigma^g \leq G$ for any $\sigma > 0$. In addition, for any
\[\sigma > \max\left\{-H_\sigma^g, 0\right\},\]
the rare event potential $V_\sigma$ verifies~\eqref{cond:log-cc-log-smooth} with $\alpha = \alpha_\sigma > 0$ and $L = L_\sigma < +\infty$, where
\[\alpha_\sigma = \frac{H^g_\sigma}{\sigma} + 1\]
and
\[L_\sigma = \left(\frac{e^{(-\mu + r)/\sigma}}{\sigma^2}K^2+\frac{2G}{\sigma} +1\right).\]
\end{prop}
\begin{rk} \label{rk:hess-g}
We observe that the condition $\sigma > \max\left\{-H_\sigma^g, 0\right\}$ implies that the smallest $\sigma$ for Proposition~\ref{prop:all-conditions} to hold depends on $g$. For instance, if $g$ of class $ C^2(\R^d, \R)$ is in addition convex, which is equivalent to $\inf_x \lambda_{\min}(\Hess(g)(x)) \geq 0$, then $H^g_\sigma \geq 0$, so any $\sigma > 0$ is admissible for Proposition~\ref{prop:all-conditions}. In contrast, if $\inf_x \lambda_{\min}(\Hess(g)(x))$ is negative, then $\sigma$ must be large enough to counteract this effect. This can be intuitively understood by observing the form of $V_\sigma$, which we recall as for $x \in \R^d$ as
$V_\sigma(x) = - \log F_\sigma(x) + \lVert x \rVert^2/2$. Note that $-\log F_\sigma(x) \to 0$ as $\sigma \to +\infty$, and $-\log F_\sigma(x) \to +\infty$ if $\sigma \to 0$. Then, for a sufficiently large $\sigma$ (a very smooth optimal density), $V_\sigma$ is merely a small perturbation of the quadratic potential $\lVert x \rVert/2$, which is a potential satisfying~\eqref{cond:log-cc-log-smooth}. If $\sigma$ is small, then $V_\sigma$ is largely contributed by $-\log F_\sigma$, so the verification of~\eqref{cond:log-cc-log-smooth} depends essentially on $g$. If $g$ is convex, then $- \log F_\sigma$ is convex for any $\sigma > 0$, so $V_\sigma$ is convex even for small $\sigma$. If $g$ is non-convex, then there is no a priori reason for $-\log F_\sigma$ to be convex, so a small $\sigma$ cannot be chosen to ensure the convexity of $V_\sigma$.
\end{rk}
Finally, we combine Proposition~\ref{prop:all-conditions} and Theorem~\ref{thm:bias-ula} in the following corollary.
\begin{corollary}\label{cor:bias-ula-sigma} 
 In the setting of Proposition~\ref{prop:all-conditions}, we fix some
    $\sigma > \max\left\{-H_\sigma^g, 0\right\}.$ If
    $\epsilon \leq 2/(\alpha_\sigma + L_\sigma)$, then for $k \geq 0$,
    \begin{equation} \label{eq:cond-timestep}
    W_2(\pi_k, \pi_*) \leq (1-\alpha_\sigma\epsilon)^k W_2(\pi_0,\pi_*) +    1.82(L_\sigma/\alpha_\sigma)(\epsilon d)^{1/2}.
    \end{equation}
\end{corollary}
The corollary shows the two roles that the smoothing parameter $\sigma$ plays in the bias of ULA. First, the maximum size of the discretization time-step $\epsilon$ in ULA is limited by the smallness of $\sigma$ to control the bias of ULA. Secondly, the bound on the bias of ULA itself is influenced by $\sigma$, but the dependence is not immediately clear. To illustrate the dependence, consider the example where $g$ is convex. We take the largest possible step size $\epsilon_\sigma = 2/(\alpha_\sigma + L_\sigma)$, which is desired in practice to speed up the convergence.
Then, Corollary~\ref{cor:bias-ula-sigma} gives for any $\sigma > 0$,
    \begin{align*}
    W_2(\rho_k, \rhosig) \leq \left(1-\frac{2}{1+L_\sigma/\alpha_\sigma}\right)^k W_2(\rho,\rhosig)  + 1.82(2d)^{1/2}\frac{L_\sigma}{\alpha_\sigma}\left(\frac{1}{\alpha_\sigma + 1} \right)^{1/2}.
    \end{align*}
From Proposition~\ref{prop:all-conditions} and the fact that $ 0 \leq H_\sigma^g \leq \sup_{x \in \R^d}\lambda_{\min}(\Hess(g)(x))$, we have under~\eqref{cond:sigma-mu}, for any $a > 0$,
\[\lim_{\sigma \to 0}\frac{L_\sigma}{(\alpha_\sigma)^a} = \lim_{\sigma \to 0}\frac{\sigma^2 + 2G \sigma +e^{(-\mu + r)/\sigma }K^2}{(\sigma^2 + H_\sigma^g \sigma)^a} = \lim_{\sigma \to 0} \frac{e^{r/\sigma}K^2}{(\sigma^2 + H_\sigma^g \sigma)^a} \to +\infty, \]
so we see that the above bias bound tends to $+\infty$ as $\sigma \to 0$, which is coherent with the fact that one cannot sample from $\rhoopt$ using ULA. The fact that $L_\sigma/(\alpha_\sigma)^2$ tends to $+\infty$ as $\sigma \to 0$ implies that $2/(\alpha_\sigma + L_\sigma)$ tends to $0$. From the requirement $\epsilon \leq 2/(\alpha_\sigma + L_\sigma)$ in Corollary~\ref{cor:bias-ula-sigma}, we see that the maximum time-step size for Corollary~\ref{cor:bias-ula-sigma} to hold decreases as $\sigma$ decreases.

On the other hand, $\alpha_\sigma, L_\sigma \to 1$ as $\sigma \to +\infty$ (and $\rhosig \to \rho$ pointwise as $\sigma \to +\infty$). The left term in the bound disappears and this gives that for any $k \geq 0$,
  \begin{align*}
    W_2(\rho_k, \rho) \leq 1.82 \ d^{1/2},
\end{align*}
which is finite but implies a bias bound that scales with the dimension $d$.
\subsection{Numerical verification for ULA} \label{subs:ULA-numerics}
\begin{figure}
    \centering
    \begin{subfigure}{0.49\textwidth}
        \includegraphics[width=\linewidth]{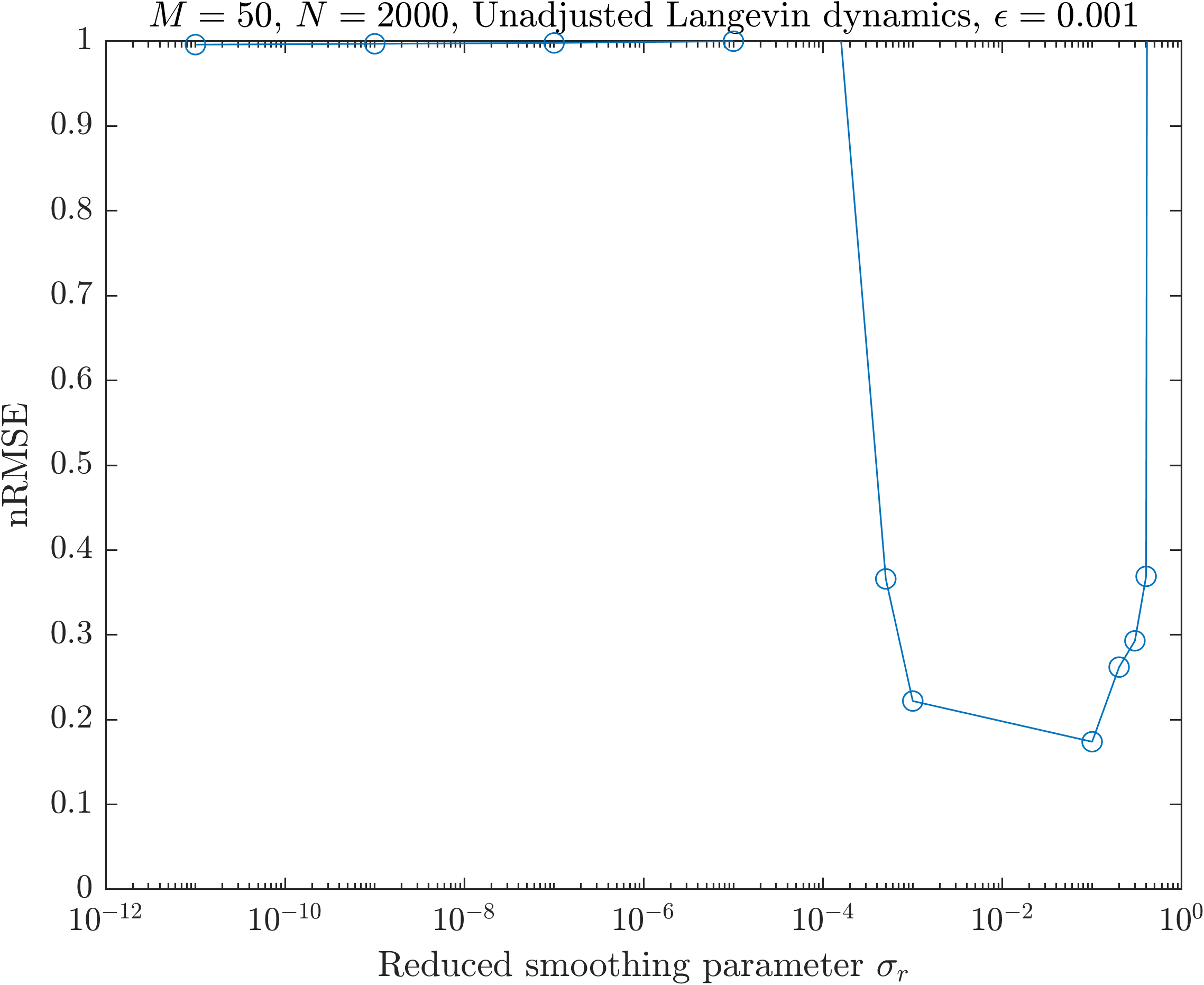}\caption{time-step size $\epsilon = 10^{-3}$}\label{fig:sigma-varies-ULA-1}
    \end{subfigure}\begin{subfigure}{0.49\textwidth}\includegraphics[width=\linewidth]{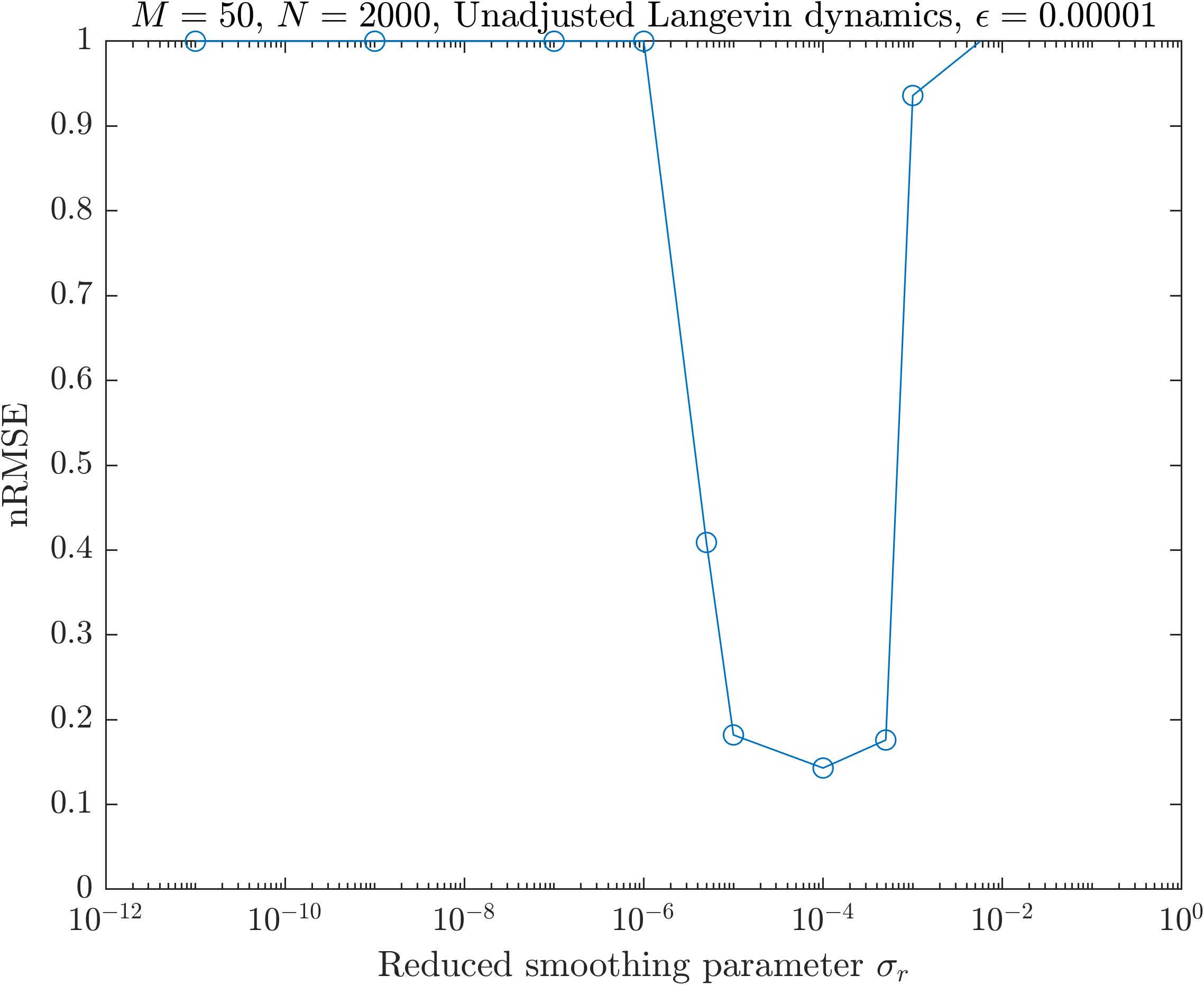}
    \caption{time-step size $\epsilon = 10^{-5}$}\label{fig:sigma-varies-ULA-2}
    \end{subfigure}
    \caption{Evolution of nRMSE with respect to the reduced smoothing parameter $\sigma_r$, for ULA with different time-step sizes $\epsilon$. Figure~\ref{fig:sigma-varies-ULA-1} shows the result for $\epsilon = 10^{-3}$ while Figure~\ref{fig:sigma-varies-ULA-2} for $\epsilon = 10^{-5}$.}
    \label{fig:sigma-varies-ULA}
\end{figure}

We investigate the trade-off phenomenon on the smoothing parameter with ULA on the classical linear LSF example, defined for all $x \in \R^d$ as
\begin{align} \label{eq:linear-LSF}
    g(x) = 5 - \frac{1}{\sqrt{d}}\sum_{i=1}^{d} x_i \text{ and } \grad g (x)= - \frac{1}{\sqrt{d}}(1,\ldots, 1)^\top.    
\end{align}
The probability $p$ to be estimated is independent of the dimension and its reference value is $p = 2.87\cdot 10^{-7}$.
We fix the dimension as $d = 100$, the number of particles as $M=50$, and the sample size for importance sampling as $N = 2000$.  We consider two time-step sizes $\epsilon = 10^{-3}$ and $\epsilon = 10^{-5}$, and for this example, we consider a deterministic stopping time which is, respectively, $T=0.1$ and $T=10^{-3}$, so the total number of time-steps is $100$. Then, the total number of LSF calls for ULA coupled with the importance sampling is $7000$. 
Inspired by the recent work in~\cite[Section 3.1]{ehre2024steinvariationalrareevent}, we study $\mu$ and $\sigma$ in the form of
\begin{equation} \label{eq:sigma-mu-forms}
 \sigma = \frac{\sqrt{3}\sigma_r}{\pi} \ \text{ and } \mu = -\log 9\left(\frac{3\sigma_r}{\pi}\right)^{1/2},    
\end{equation}
with $\sigma_r$ ranging from $10^{-11}$ to $1$.
The results are plotted in Figure~\ref{fig:sigma-varies-ULA}.
For both time-step sizes, we see that when $\sigma_r$ is too small or too large, the nRMSE goes to $1$. When the time-step size is $10^{-5}$ in Figure~\ref{fig:sigma-varies-ULA-2}, the admissible range of the smoothing parameter outside of which the nRMSE goes to $1$, shifts to a smaller range of values compared to that of $\epsilon = 10^{-3}$. This is in line with Remark~\ref{rk:hess-g}, where if the smoothing parameter is too small, then the time-step size could not be too high to ensure the condition $\sigma > \max\{-H^g_\sigma , 0 \}$ holds. For the larger time-step size $\epsilon = 10^{-3}$ the optimal value among those tested is $\sigma_r = 0.1$ and the admissible range of parameters appears to be $[5\cdot 10^{-4}, 0.4]$, while for the smaller time-step size $\epsilon = 10^{-5}$, the optimal value is $10^{-4}$ and the admissible range is $[5\cdot 10^{-6}, 5\cdot 10^{-3}]$. This is coherent with the condition $\epsilon \leq 2/(\alpha_\sigma + L_\sigma)$ which indicates that a smaller $\sigma$ requires a smaller maximum time-step size (because $L_\sigma/\alpha_\sigma \to +\infty$ as $\sigma \to 0$). 

\subsection{Stabilization w.r.t.\ the smoothing parameter via adaptive time-step size and stopping criterion}
We want to check if the adaptive time-step size and stopping criterion presented in Section~\ref{sec:adaptive-timestep} and~\ref{sec:stopping-criterion} stabilize the nRMSE w.r.t.\ the smoothing paramter. 
We present a preliminary illustration on ULA that our proposed adaptive time-step size and stopping criterion leads to enhanced robustness with respect to the smoothing parameter choice.

\begin{figure}
    \centering
    \includegraphics[width=0.5\linewidth]{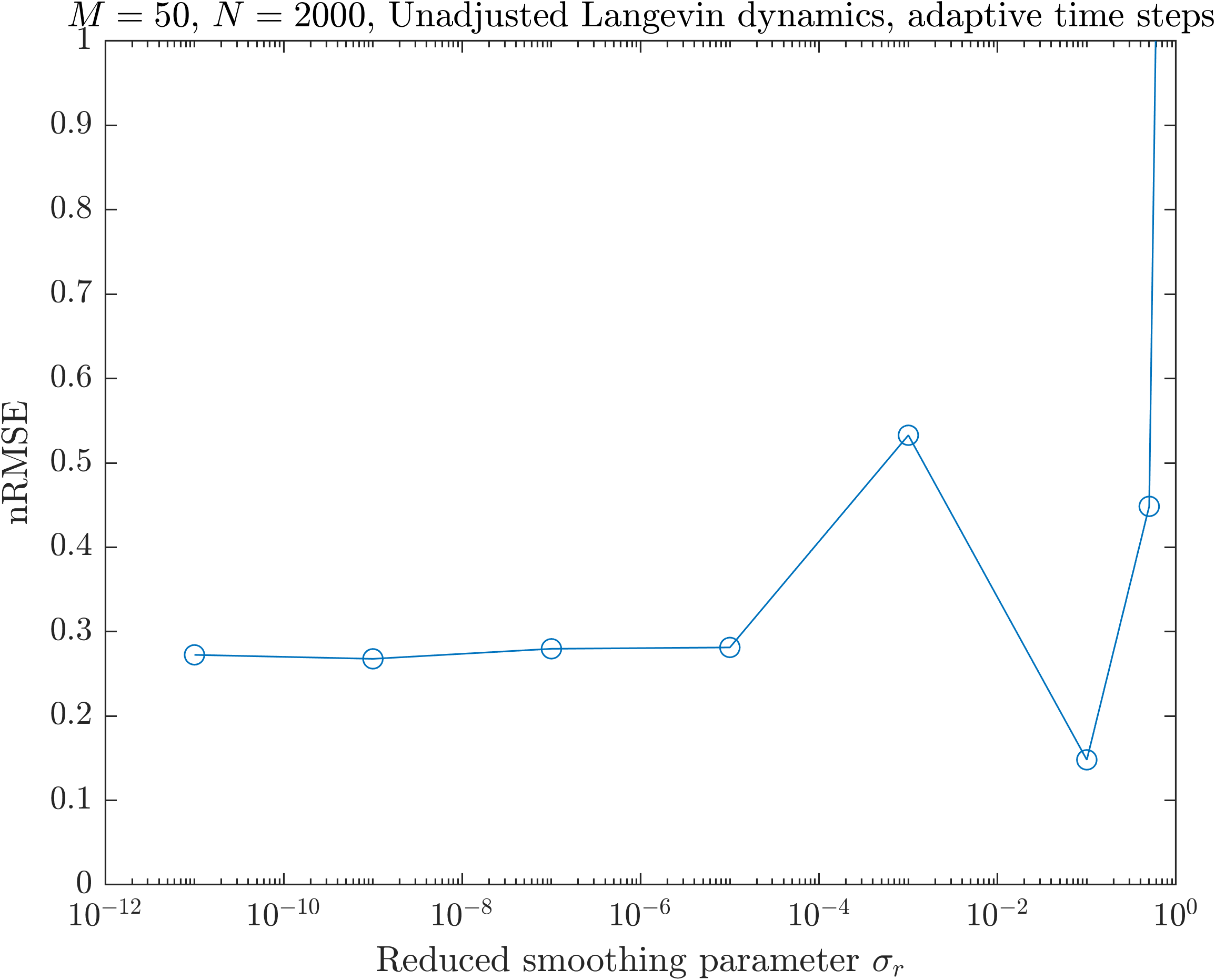}
    \caption{Evolution of nRMSE with respect to the reduced smoothing parameter $\sigma_r$, for ULA with adaptive time-step size and stopping criterion described in Section~\ref{sec:adaptive-timestep} and~\ref{sec:stopping-criterion}.}
    \label{fig:sigma-varies-ULA-adaptive}
\end{figure}
In Figure~\ref{fig:sigma-varies-ULA-adaptive}, we observe that the nRMSE stabilizes at $0.3$ for very small smoothing parameters starting from $10^{-5}$. On the other hand, the adaptive time-step size and stopping criterion do not protect against the effect of too large choices of $\sigma$. 

We present a brief study on the influence of $\sigma$ on the nRMSE of the estimator to check if the phenomenon of the trade-off with respect to the smoothing parameter for ULA, translates into ALDI-IS with DBSCAN with the adaptive criteria. We study again $\mu$ and $\sigma$ in the form of~\eqref{eq:sigma-mu-forms},
with this time $\sigma_r \in \{10^{-11}, 10^{-9}, 10^{-7}, 10^{-5}, 10^{-3}, 10^{-1}, 0.5, 1\}$.
We consider again the hyperplane example, and an additional four branches examples, defined for any $x \in \R^2$,
\begin{align} \label{eq:fourbranches-LSF}
g(x) = \min \begin{pmatrix}
        0.1(x_1-x_2)^2 - (x_1 + x_2)/\sqrt{2} + 3, \\
        0.1(x_1- x_2)^2 + (x_1 + x_2)/\sqrt{2} + 3, \\
        x_1 - x_2 + 7/\sqrt{2}, \\
        x_2 - x_1 + 7 /\sqrt{2},
    \end{pmatrix} := \min\left(g_1(x),g_2(x), g_3(x), g_4(x)\right).    
\end{align}

\begin{figure}
    \centering
    \begin{subfigure}{0.49\textwidth}
        \includegraphics[width=\linewidth]{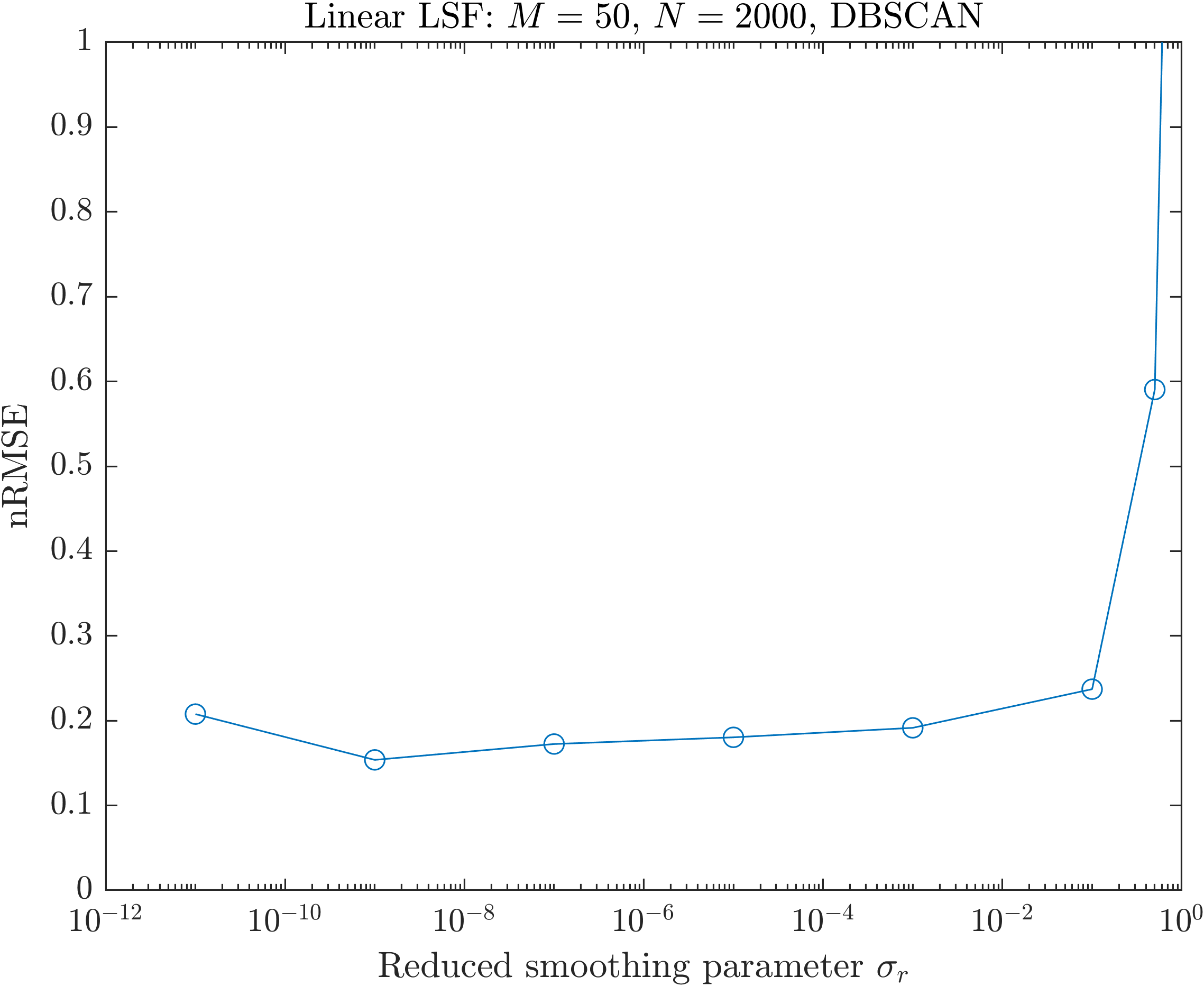}
        \caption{Linear LSF}
    \end{subfigure}
    \begin{subfigure}{0.49\textwidth}\includegraphics[width=\linewidth]{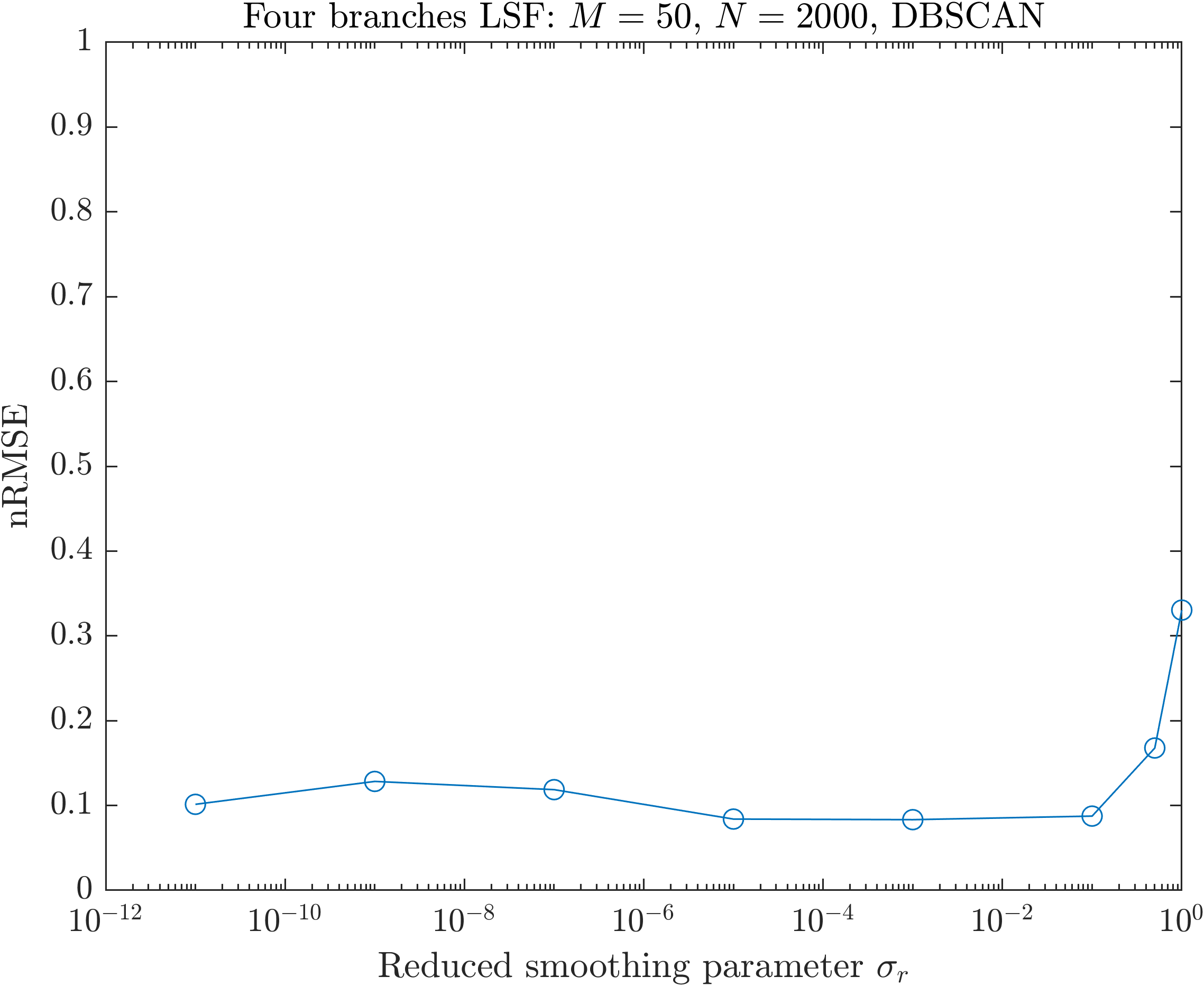}
    \caption{Four branches LSF}
    \end{subfigure}
    \caption{Evolution of nRMSE with respect to the reduced smoothing parameter $\sigma_r$, for linear LSF and four branches LSF.}
    \label{fig:sigma-varies}
\end{figure}

The reader is referred to Section~\ref{sec:four-branches} for further details of implementation for the four branches example. The results are plotted in Figure~\ref{fig:sigma-varies}.
\begin{figure}
    \centering
    \begin{subfigure}{0.49\textwidth}
        \includegraphics[width=\linewidth]{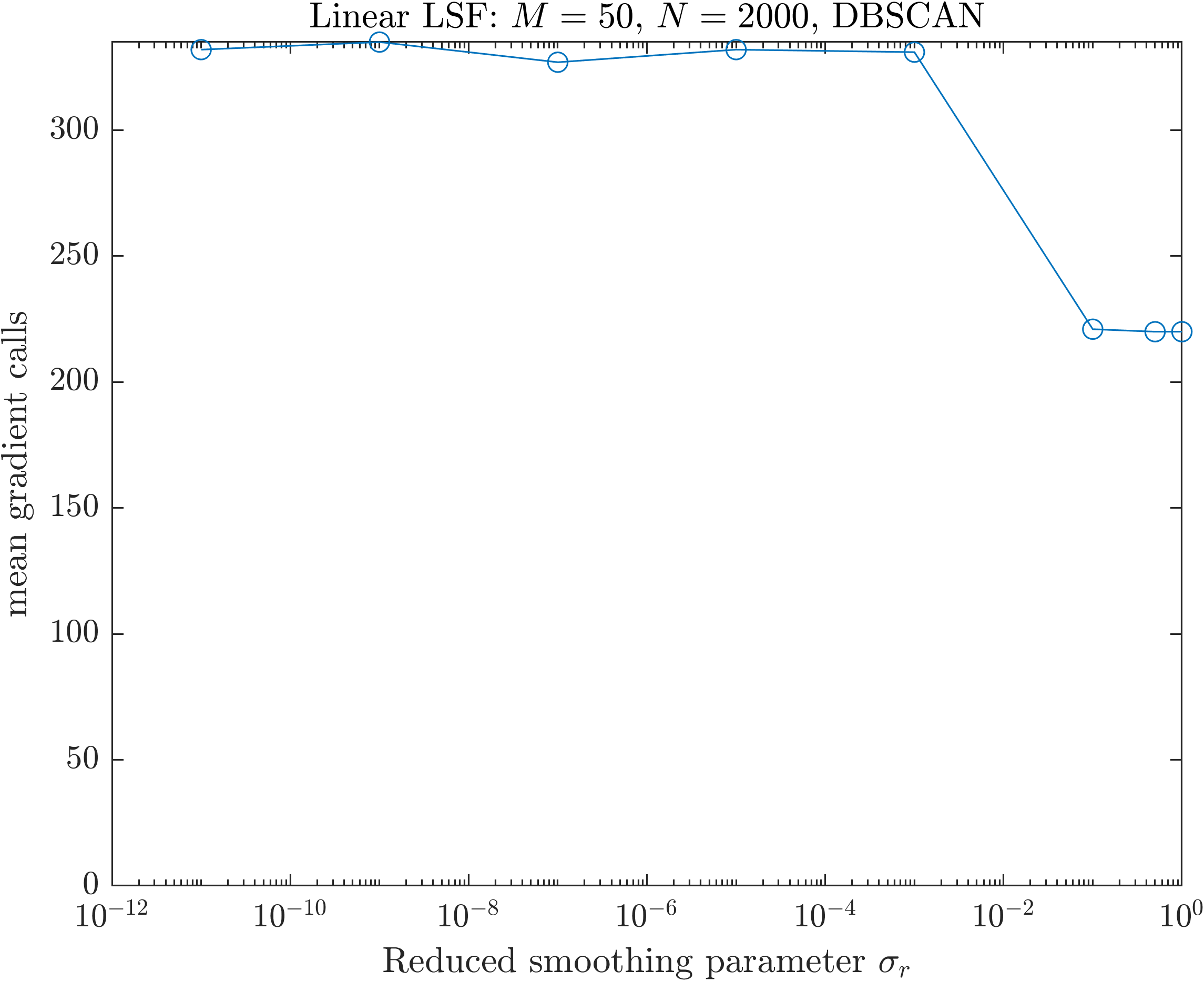}
        \caption{Linear LSF}
    \end{subfigure}
    \begin{subfigure}{0.49\textwidth}\includegraphics[width=\linewidth]{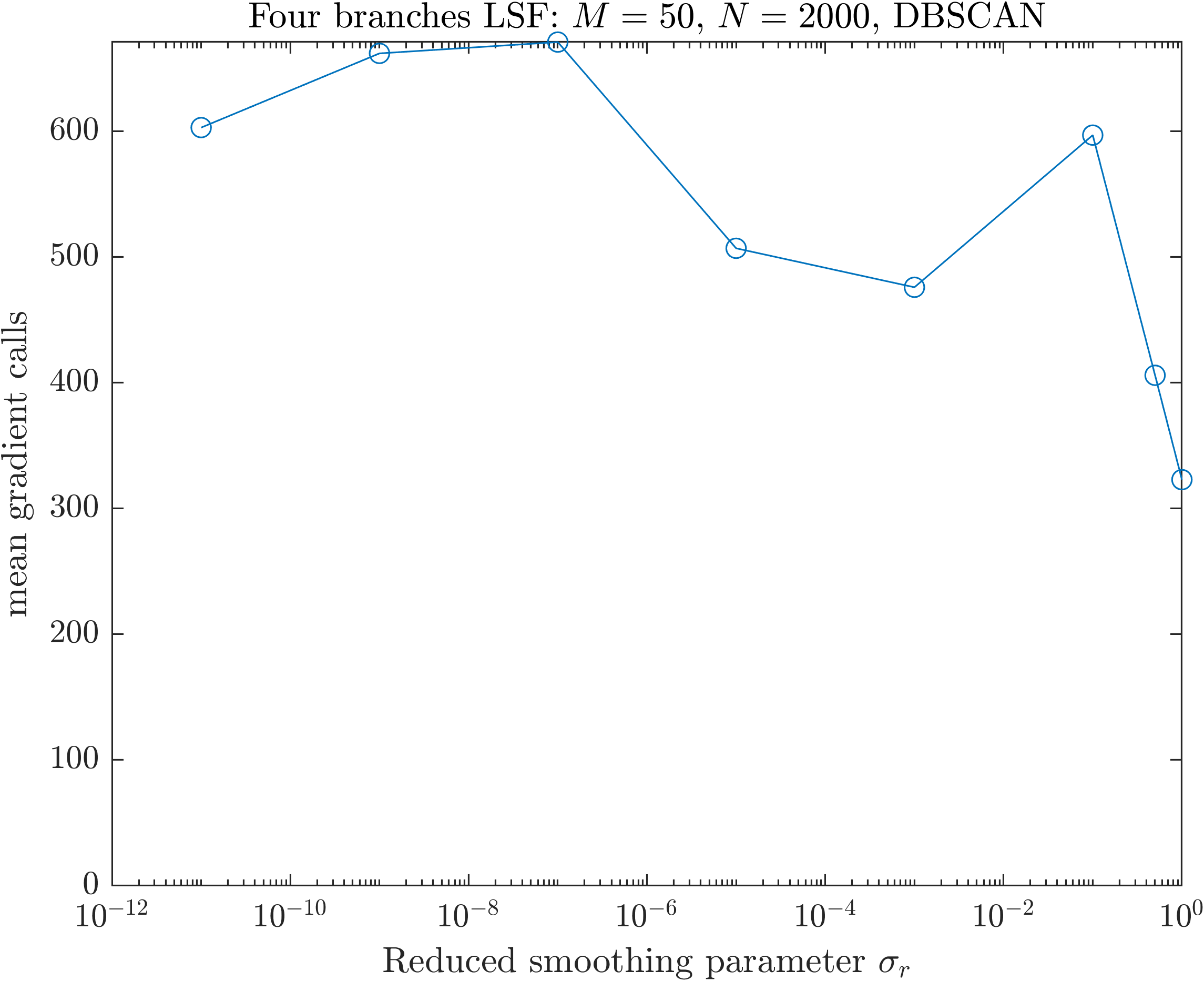}
    \caption{Four branches LSF}
    \end{subfigure}
    \caption{Evolution of the mean gradient calls with respect to the reduced smoothing parameter $\sigma_r$, for linear LSF and four branches LSF.}
    \label{fig:sigma-varies-gradientcalls}
\end{figure}
The difference between these two examples is that the linear LSF is $1$-Lipschitz continuous on $\R^d$, since for any $x,y \in \R^d$,
\[ \lvert g(y) - g(x) \rvert = \frac{1}{\sqrt{d}}\left\lvert\sum_{i=1}^d (y_i - x_i) \right\rvert \leq \frac{1}{\sqrt{d}}\sum_{i=1}^d \left\lvert y_i - x_i \right\rvert \leq \left\lVert y - x \right\rVert, \]
where we have used the fact that the $1$-norm on $\R^d$ verifies for any $z \in \R^d$
\[\sum_{i=1}^d \lvert z_i \rvert \leq \sqrt{d}\lVert z \rVert.\]
On the contrary, the four branches LSF is not Lipschitz continuous on $\R^d$, since it does not verify the property of being differentiable almost everywhere. 
Remark~\ref{rk:hess-g} suggests that when the LSF $g$ is Lipschitz continuous on $\R^d$, the optimal smoothing parameter would be smaller compared to another LSF that is not. We can observe this effect in Figure~\ref{fig:sigma-varies}, in which among the $\sigma_r$'s tested, the optimal $\sigma_r$ is $10^{-9}$ for the linear LSF, while it is $10^{-3}$ for the four-branches LSF. The phenomenon of the trade-off on $\sigma$ between the importance sampling error and the ease of sampling is lightly present for ALDI-IS with DBSCAN, albeit it is clear that a smaller than optimal $\sigma$ does not degrade the nRMSE as much as a $\sigma$ which is too large. Compared to ULA in Figure~\ref{fig:sigma-varies-ULA}, ALDI-IS with DBSCAN is more robust with respect to the effects of small $\sigma$. This is largely due to the adaptive time-step scheme we chose in~\eqref{eq:adaptive-time-step}. If the smoothing parameter is too small, then the size of time-steps is adjusted to be small in order to preserve the stability of the scheme. In fact, the adaptive size of the time-steps coupled with our stopping criterion preserve a relatively stable number of gradient calls for ALDI, therefore making the nRMSE robust w.r.t.\ the smoothing parameter. The stable number of gradient calls is illustrated in Figure~\ref{fig:sigma-varies-gradientcalls}.

\section{Numerical experiments} \label{sec:num}
In this section, we illustrate our proposed algorithm ALDI-IS (with and without DBSCAN) on three standard test cases found in the literature. We compare our results with other algorithms which will be introduced in their respective test cases. For the expectation-maximization algorithm for fitting vMFNM model, as well as the improved cross-entropy with vMFNM model from~\cite{papaioannou_improved_2019}, we use the free software Reliability Analysis Tools from the Engineering Risk Analysis Group of the Technical University of Munich, available at \url{https://github.com/ERA-Software/Overview}.

\subsection{Linear LSF with dimension scaling} \label{sec:linear}
\begin{figure}
    \centering
    \begin{subfigure}{0.49\textwidth}
        \includegraphics[width=\linewidth]{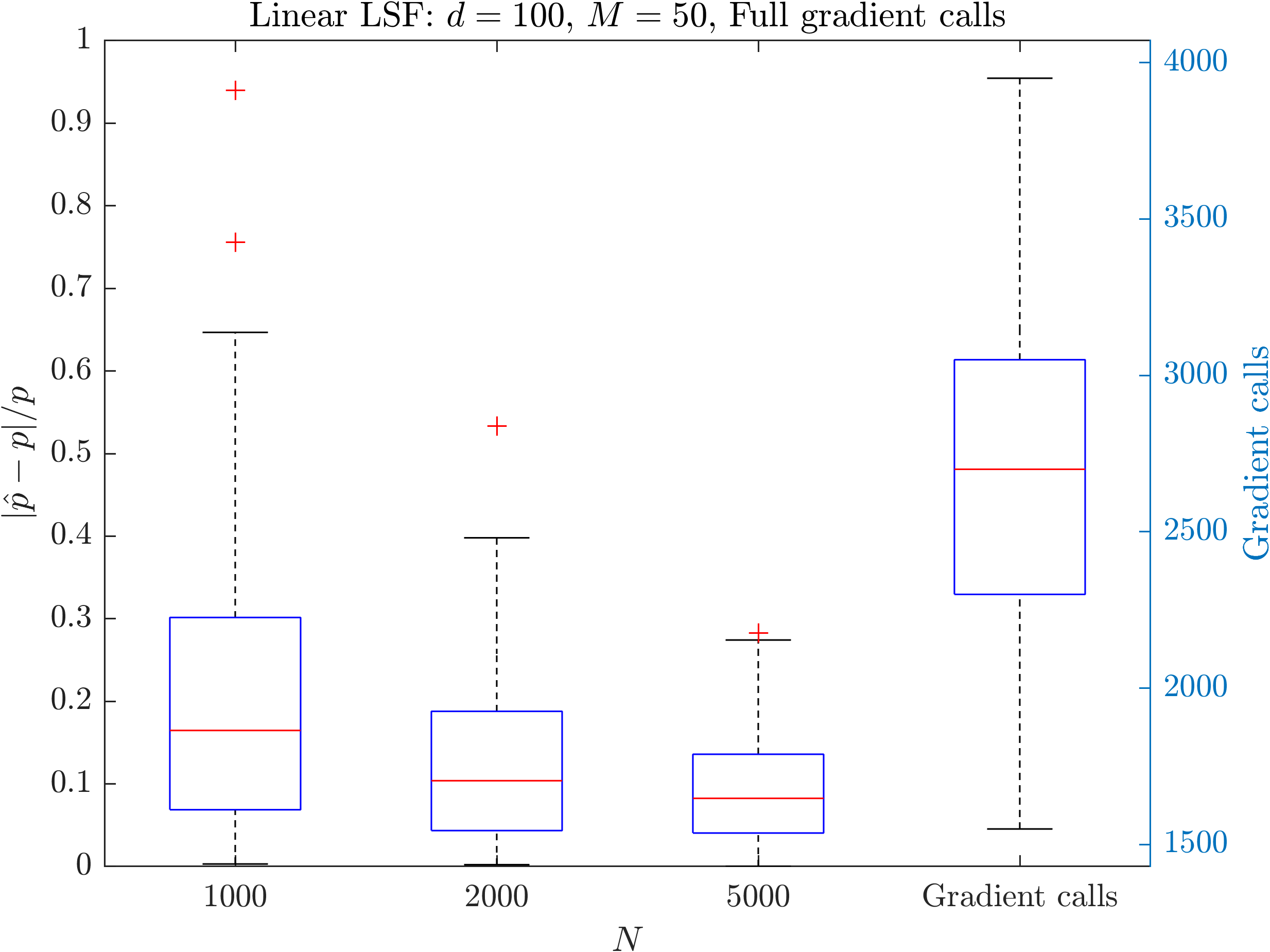}
    \caption{Full gradient calls}\label{fig:linear_LSF-1}    
    \end{subfigure}
    \begin{subfigure}{0.49\textwidth}\includegraphics[width=\linewidth]{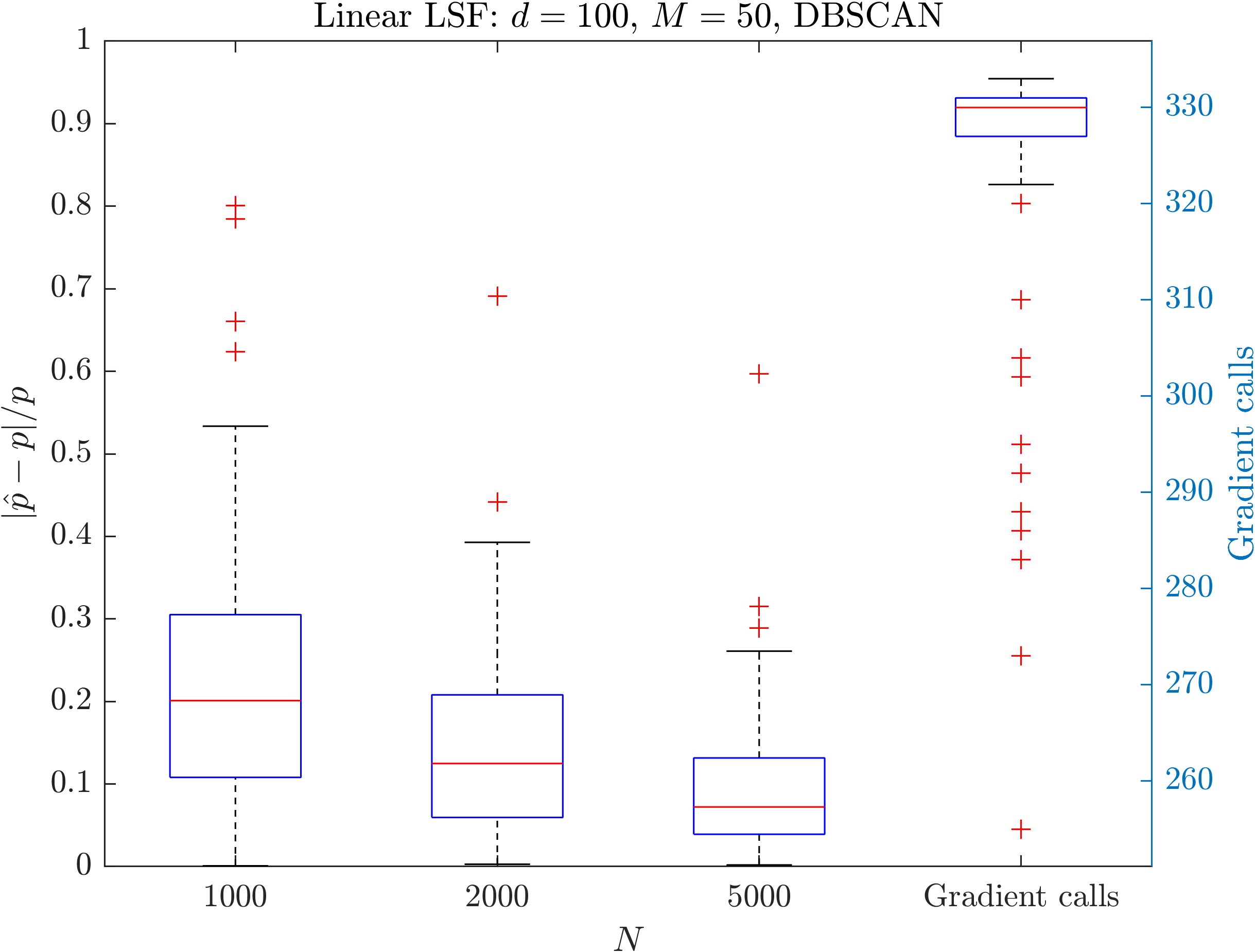} 
        \caption{DBSCAN}\label{fig:linear_LSF-2}  
    \end{subfigure}
    \caption{Boxplots of the normalized absolute error as a function of the number of IS samples $N$ for linear LSF. Figure~\ref{fig:linear_LSF-1} uses ALDI-IS without DBSCAN (full gradient calls), while Figure~\ref{fig:linear_LSF-2} uses ALDI-IS with DBSCAN.}
    \label{fig:linear_LSF}
\end{figure}

First, we consider again the classical linear LSF example which was defined in~\eqref{eq:linear-LSF}.
The probability $p$ to be estimated is independent of the dimension and its reference value is $p = 2.87\cdot 10^{-7}$.
We fix the dimension as $d = 100$, and the following parameters for ALDI:
$Q=4$ (we recall that this means that we approach $A$ sequentially with 4 ALDI's, each using the final particles from the previous dynamics as new initial particles),
$(q_i) = (1,0.5,0.05,0)$, $(\gamma_i) = (1,0.5,0.01,10^{-3})$, $(\eps^{\cumu}_i) =(0.1,0.1,0.1,0.01)$ and $M=50$.
Once the samples are collected, they are fitted with the expectation-maximization algorithm to obtain $\rho_{\hat{\boldsymbol{\nu}}}^{\vMFNM}$, from which we generate $N= 1000, 2000$ and $5000$ i.i.d.\ samples for the importance sampling step. We denote the estimated probability by $\hat p$ for simplicity. The result is plotted as a boxplot with $100$ repetitions of the whole ALDI-IS scheme in Figure~\ref{fig:linear_LSF-1}. The $x$-axis represents the number of importance sampling samples $N$, and the left $y$-axis represents the relative error $\lvert \hat p - p \rvert/p$. The right $y$-axis represents the total number of gradient calls of $g$ used in the ALDI step. We note that the number of calls of $g$ itself is the sum of the importance sampling sample size $N$ 
plus the gradient calls of $g$ during ALDI, since one evaluation of the gradient of $V_\sigma$ 
requires one evaluation of $g$ and $\grad g$ according to~\eqref{eq:grad-V}.

In addition, 
we emphasize that in this example, for all $x \in \R^d$ the value $\grad g(x)$ does not depend on $x$. Therefore, ALDI could be run with only one gradient evaluation per iteration instead of $M$ gradient evaluations. In practical scenarios, such a specific property of $\grad g$ cannot be known beforehand, but we expect ALDI-IS coupled with DBSCAN to work well in this test case. The results are plotted in Figure~\ref{fig:linear_LSF-2}, in which we can observe that the number of gradient calls is drastically decreased. The mean number of gradient calls is decreased by an order of magnitude, from around 2500 (Figure~\ref{fig:linear_LSF-1}) to around 300 (Figure~\ref{fig:linear_LSF-2}). Note that the number of gradient calls is not decreased by a factor of $50$ since the stopping time $T$ is random and is realized later with DBSCAN. This test case represents the ideal scenario for the clustering algorithm to work well. 

For this test case, we compare the nRMSE of our estimator with readily available results in the literature, including component-wise Metropolis-Hastings subset simulation (MH-SuS)~\cite{PAPAKONSTANTINOU2023}, rejection-based Hamiltonian Monte Carlo subset simulation (HMC-SuS)~\cite{WANG2019} and Hamiltonian Markov chain Monte Carlo as approximate sampling target with post-processing adjustment (HMCMC-ASTPA)~\cite{PAPAKONSTANTINOU2023}. The results are plotted in Table~\ref{table:linear}.
\begin{table}
\centering
\begin{tabular}{|c|c|c|c|}
 \hline
 Method & nRMSE & LSF calls & LSF gradient calls 
 \\ \hline
 MH-SuS~\cite{PAPAKONSTANTINOU2023}  & 0.69 & 6443& 0 
 \\HMC-SuS~\cite{WANG2019}   & 0.28 &6403 & 6403
 \\HMCMC-ASTPA~\cite{PAPAKONSTANTINOU2023}  & 0.12  &2225& 2225
 \\ALDI-IS ($N=2000, M=50$) & 0.19&4680& 2680
 \\ALDI-IS DBSCAN ($N=2000, M=50$)& 0.18 &2325& 325
 \\ALDI-IS DBSCAN ($N=2000, M=100$) &0.09&2525& 525
 \\
\hline
\end{tabular}
\caption{nRMSE, mean LSF calls and mean LSF gradient calls for different methods on the linear LSF with $d=100$. Results for MH-SuS, HMC-SuS and HMCMC-ASTPA are directly taken from the corresponding references.}
\label{table:linear}
\end{table}
We observe that ALDI-IS works better than the subset simulation family of algorithms. However, without DBSCAN, it has slightly higher nRMSE with higher budget than HMCMC-ASTPA. ALDI-IS with DBSCAN significantly reduces the LSF gradient calls and outperforms HMCMC-ASTPA.
As mentioned, DBSCAN works specifically well in this example since the gradient of the LSF does not depend on $x$, which allows the use of a higher number of particles ($M=d=100$) without a significant increase in LSF gradient calls. This allows to achieve a better performance than HMCMC-ASTPA. We also highlight that any new advances on clustering algorithms will benefit ALDI-IS. 

\subsection{Four branches LSF} \label{sec:four-branches}
\begin{figure}
    \centering
    \begin{subfigure}{0.3\textwidth}
        \includegraphics[width=\linewidth]{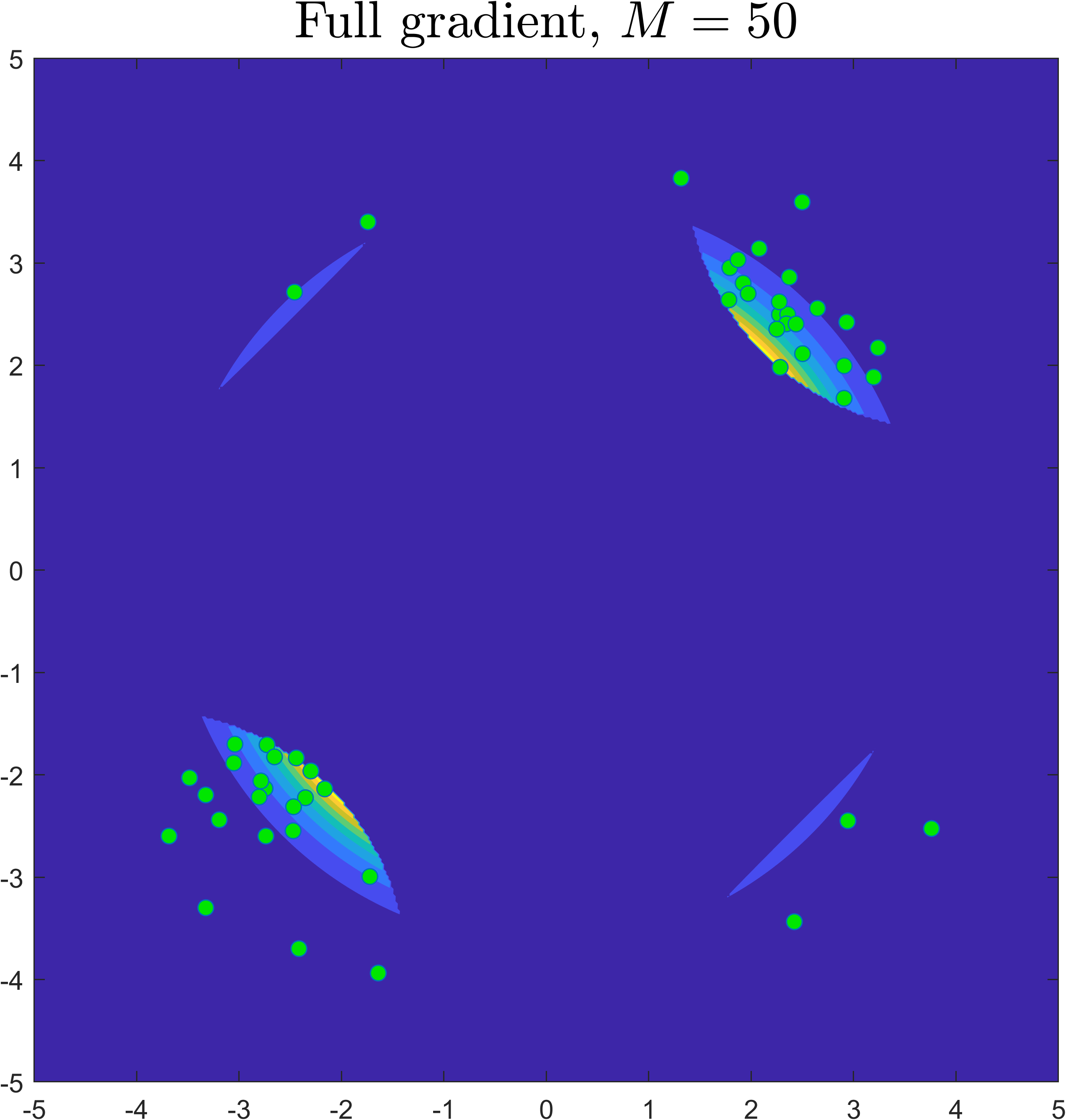}
        \caption{Full gradient calls}\label{fig:Fourbranches-1}
    \end{subfigure}
    \begin{subfigure}{0.3\textwidth}\includegraphics[width=\linewidth]{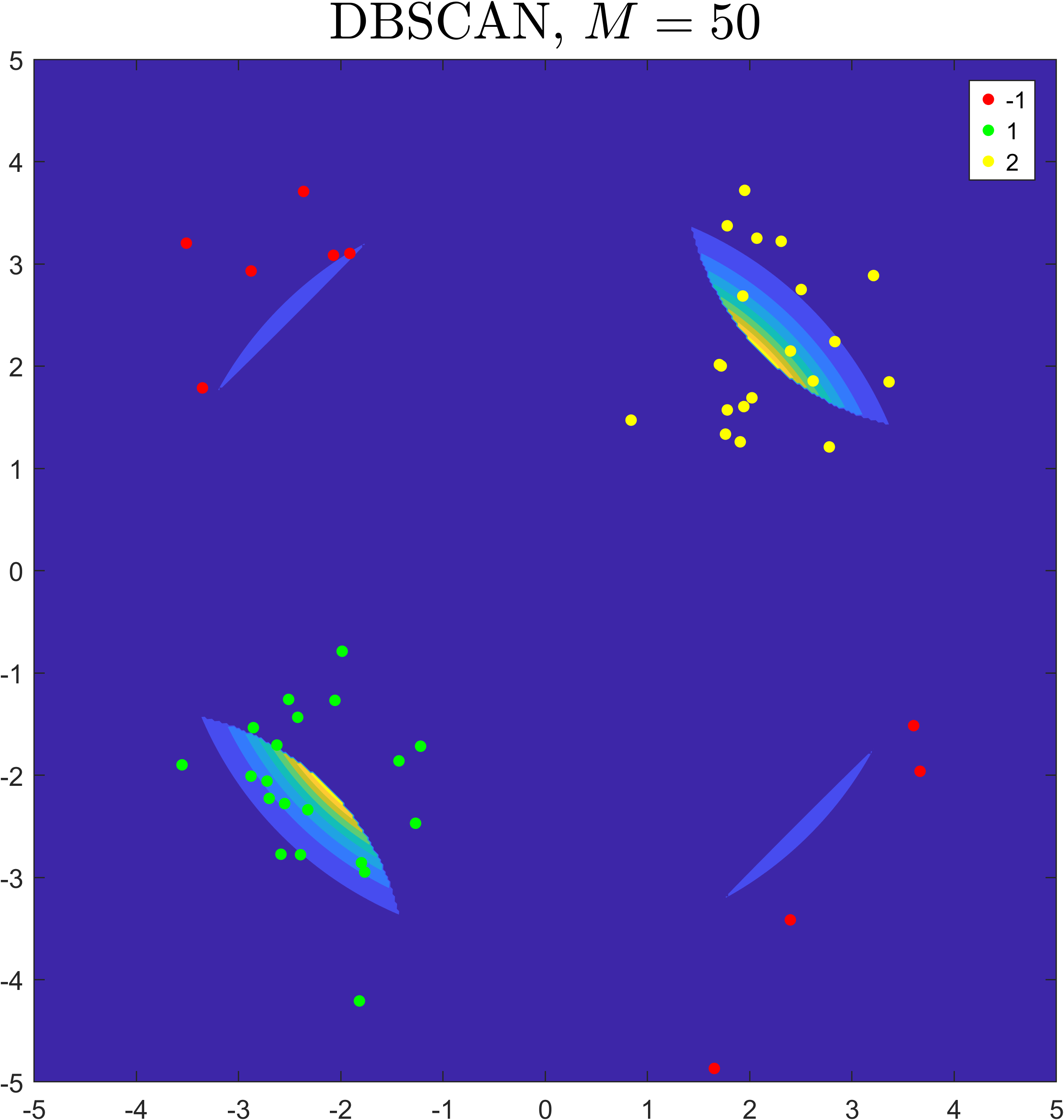}
        \caption{DBSCAN}\label{fig:Fourbranches-2}
    \end{subfigure}
    \begin{subfigure}{0.3\textwidth}\includegraphics[width=\linewidth]{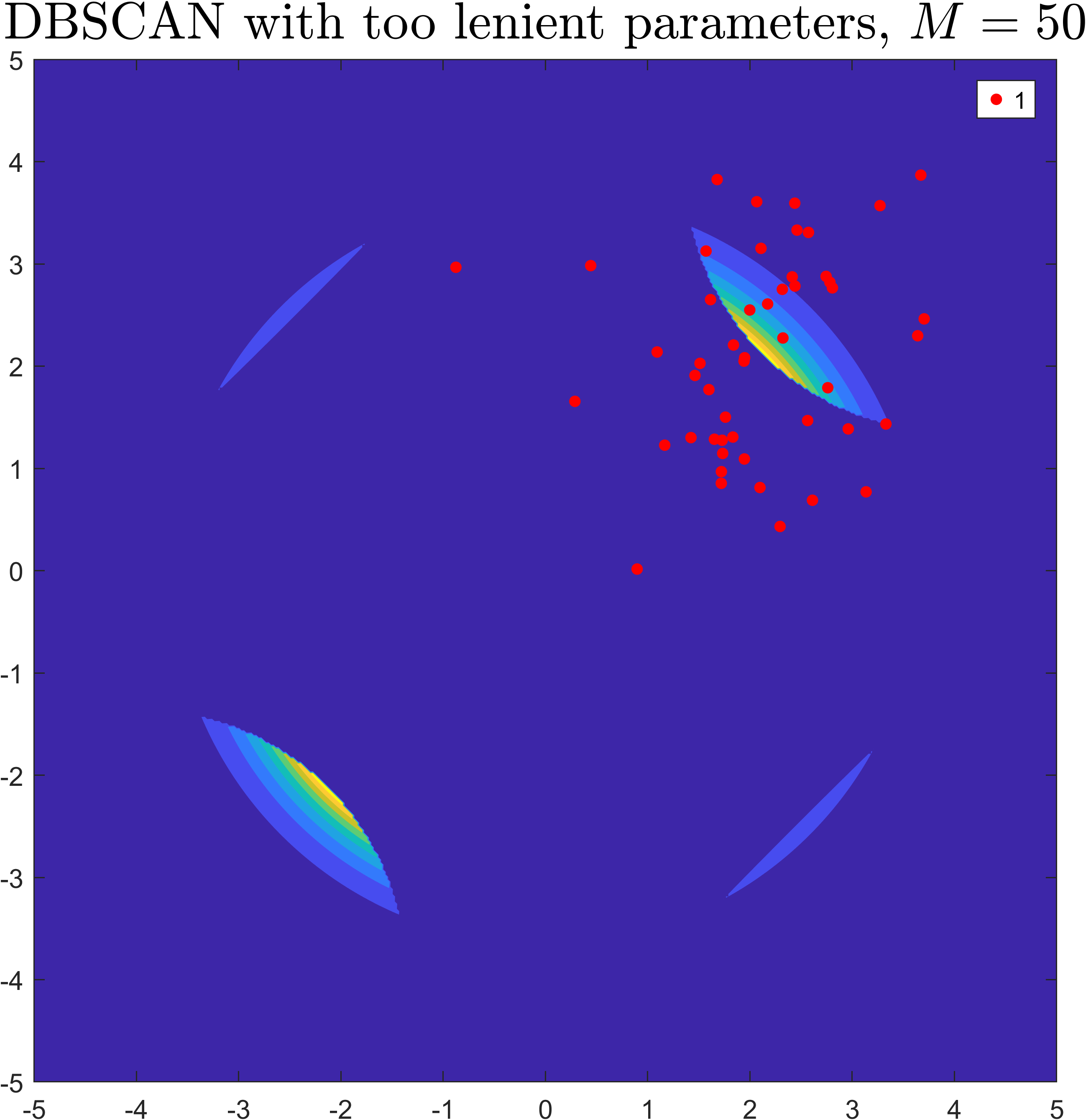}
    \caption{Too lenient DBSCAN}
        \label{fig:Fourbranches-3}
    \end{subfigure}
    \caption{Illustration of the optimal IS distribution, on which the repartition of $50$ particles in a repetition of ALDI-IS and ALDI-IS with DBSCAN are shown respectively in Figure~\ref{fig:Fourbranches-1} and~\ref{fig:Fourbranches-2}. A run for which DBSCAN had too lenient parameters $(\varepsilon = 2)$ is also plotted in Figure~\ref{fig:Fourbranches-3}.}
    \label{fig:Fourbranches}
\end{figure}

We consider the four branches LSF example, defined in~\eqref{eq:fourbranches-LSF}. 
The reference probability is $p = 2.22\cdot 10^{-3}$~\cite{papaioannou_improved_2019}. Note that the gradient of $g$ is not properly defined, so the numerical implementation uses the subgradient of $g$, still denoted as $\grad g$, given for all $x \in \R^2$ as
\[\grad g(x) = \begin{cases}
    [0.2(x_1 - x_2) - 1/\sqrt{2}, -0.2(x_1 - x_2) - 1/\sqrt{2}]^\top &\text{ if } \arg\min_{k}(g_k(x)) = 1
    \\ [0.2(x_1 - x_2) + 1/\sqrt{2}, -0.2(x_1 - x_2) + 1/\sqrt{2}]^\top &\text{ if } \arg\min_{k}(g_k(x)) = 2
    \\ [1, -1]^\top &\text{ if } \arg\min_{k}(g_k(x)) = 3
    \\ [-1, 1]^\top &\text{ if } \arg\min_{k}(g_k(x)) = 4.
    
\end{cases}\]
The optimal IS distribution is shown in the three plots in Figure~\ref{fig:Fourbranches} in the background. Figure~\ref{fig:Fourbranches-1} and~\ref{fig:Fourbranches-2} show the particles of one repetition of a complete run of ALDI-IS and ALDI-DBSCAN respectively. In the latter, the colours of the particles separate the clusters identified by DBSCAN. Red particles labelled as $-1$ represent the outlier particles identified by DBSCAN that do not belong to any cluster, so they are evolved using their own gradient instead of that of a mean particle. The parameters used for ALDI-IS are as follows: $(q_i) = (1,0.5,0.05,0)$, $(\gamma_i) = (1,0.5,0.01,10^{-3})$, $(\eps^{\cumu}_i) =(0.1,0.1,0.1,0.05)$ and $M=50$. 

In Figure~\ref{fig:Fourbranches-2}, the DBSCAN is tuned with epsilon neighbourhood $\varepsilon =1 = d/2$ as recommended by the heuristics, and minimum number of neighbours $\delta = 5$. In this case DBSCAN successfully clusters the particles on the top right and bottom left modes, and identifies the particles in the remaining two modes (with lower probability mass) as outliers since the particles are few. Regardless, this still allows to save LSF gradient calls on the two identified clusters. One can also visually observe that the particles in ALDI-DBSCAN are more dispersed in each of the modes, whereas for the full gradient call version of ALDI, all the particles are in the failure set $A$. Figure~\ref{fig:Fourbranches-3} however shows a situation when DBSCAN is tuned with too lenient parameters, in this case the epsilon neighbourhood $\varepsilon = 2 = d$. All the particles are identified as belonging in the same cluster, and thus share the same gradient which pushes the particles onto the same mode of the optimal IS distribution.

\begin{figure}
    \centering
    \begin{subfigure}{0.49\textwidth}\includegraphics[width=\linewidth]{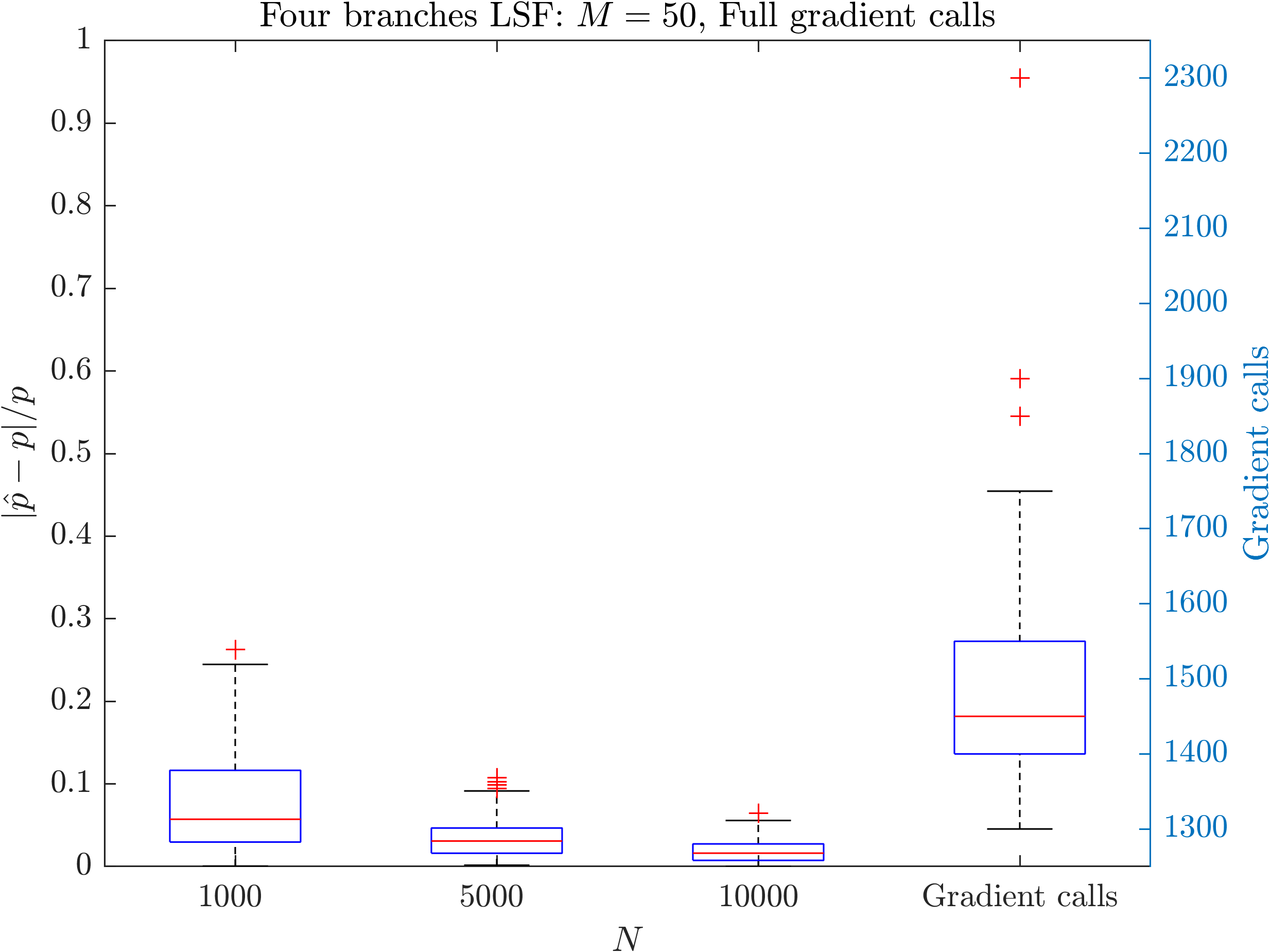}
    \caption{Full gradient calls}
    \end{subfigure}
    \begin{subfigure}{0.49\textwidth}\includegraphics[width=\linewidth]{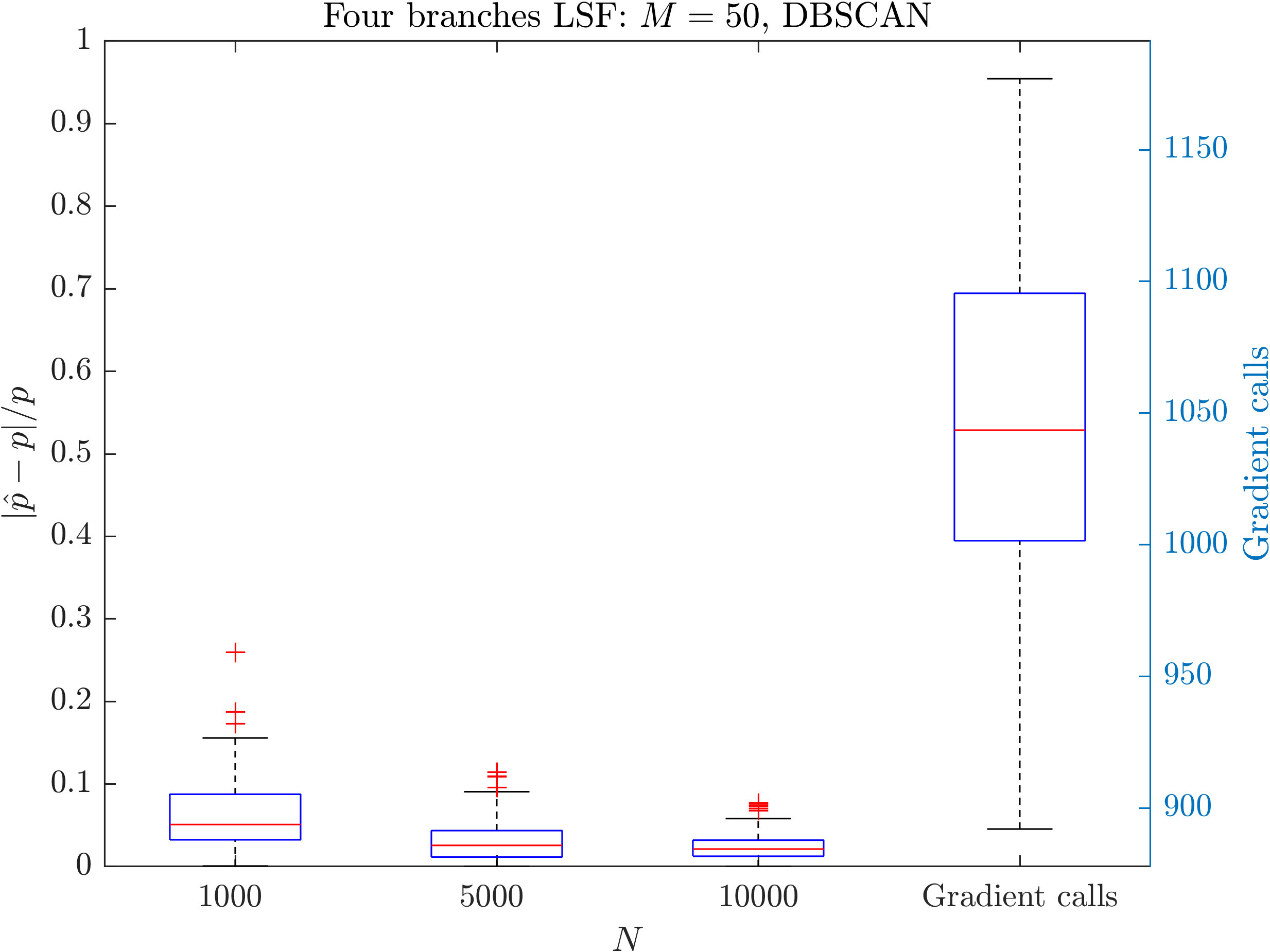}\caption{DBSCAN}
    \end{subfigure}
    \caption{Boxplots of the normalized absolute error as a function of the number of IS samples $N$ and gradient calls for four branches LSF, with and $M=50$ particles for ALDI-IS and ALDI-IS with DBSCAN respectively.}
    \label{fig:fourbranches-full}
\end{figure}

We show the error distribution and the gradient calls boxplot for ALDI-IS and ALDI-IS with DBSCAN in Figure~\ref{fig:fourbranches-full}. The existing results in the literature include~\cite[Section 6.2]{papaioannou_improved_2019} for some results on cross-entropy and improved cross-entropy scheme. However, the total number of LSF calls was not reported in the paper. Thus, we rerun improved cross-entropy with von-Mises Fisher Nakagami mixture (iCE-vMFNM) with the initial number of mixtures set to $4$, as it was the algorithm with the best result in~\cite{papaioannou_improved_2019}. We use the version in the Reliability Analysis Tools mentioned previously, 
to avoid any potential difference in implementation. The results are collected in Table~\ref{table:four-branches}. We observe that both ALDI-IS and ALDI-IS with DBSCAN outperform iCE-vMFNM with less LSF calls, although they require the gradient calls. 
We also observe that for $M=10$, ALDI-IS with and without DBSCAN provide similar results, since DBSCAN does not work properly for too few particles.

Lastly, ALDI 
automatically identifies modes in the target probability distribution, in contrast to improved cross-entropy methods with mixture model whose performance is dependent on the initial number of mixtures. The user has to specify an educated guess of the number of modes of the target probability distribution. This constituted one of our motivations to explore ALDI for rare event estimation.
\begin{table}
\centering
\begin{tabular}{|c|c|c|c|c}
 \hline
Parameters &nRMSE & LSF calls & LSF gradient calls 
 \\ \hline
 \multicolumn{4}{|c|}{iCE-vMFNM (4 initial number of mixtures)}
 \\ \hline $N_{\text{it}} = 100$  & 0.59 &4533 & 0 
 \\ $N_{\text{it}} = 250$  & 0.21 &10330 & 0 
 \\$N_{\text{it}} = 500$  & 0.17 &19650 & 0 
  \\ \hline
 \multicolumn{4}{|c|}{ALDI-IS without DBSCAN}
 \\ \hline $N=1000, M=10$ & 0.12 &1200& 200
 \\$N=1000, M=50$ &0.10&2500& 1500 
 \\ \hline
 \multicolumn{4}{|c|}{ALDI-IS with DBSCAN }
 \\ \hline $N=1000, M=10$ & 0.12 &1210& 210
 \\$N=1000, M=50$ &0.08 &2032& 1032 
\\ \hline
\end{tabular}
\caption{The nRMSE, mean LSF calls and mean LSF gradient calls for iCE-vMFNM and ALDI-IS DBSCAN on the four branches LSF, over $100$ runs.}
\label{table:four-branches}
\end{table}

\subsection{PDE example: LSF involving the discretization of a random field} \label{sec:PDE-example}
\begin{figure}
    \centering
    \begin{subfigure}{0.49\textwidth}\includegraphics[width=\linewidth]{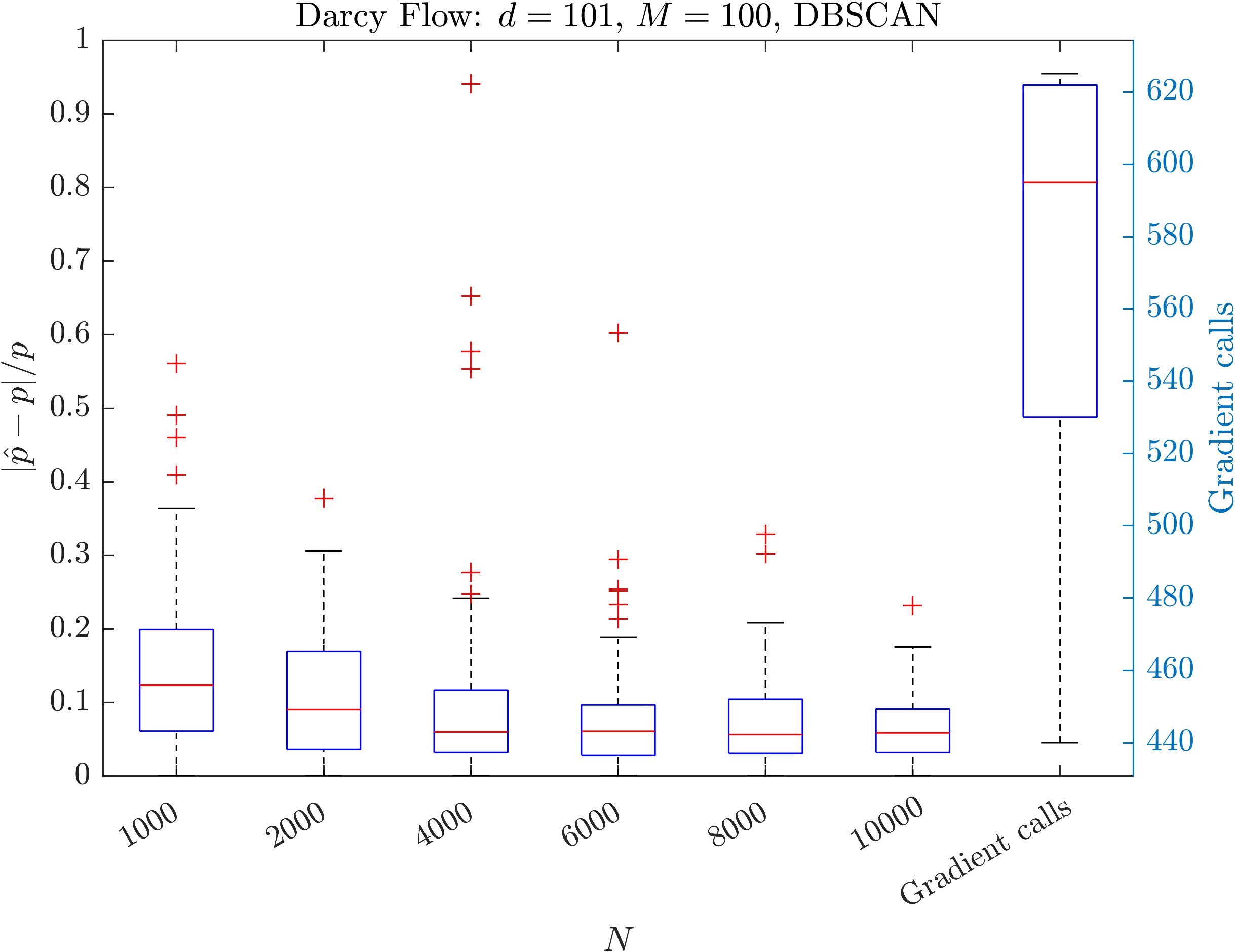}
        \caption{$M=100$} \label{fig:Darcy-1}
    \end{subfigure}
    \begin{subfigure}{0.49\textwidth}\includegraphics[width=\linewidth]{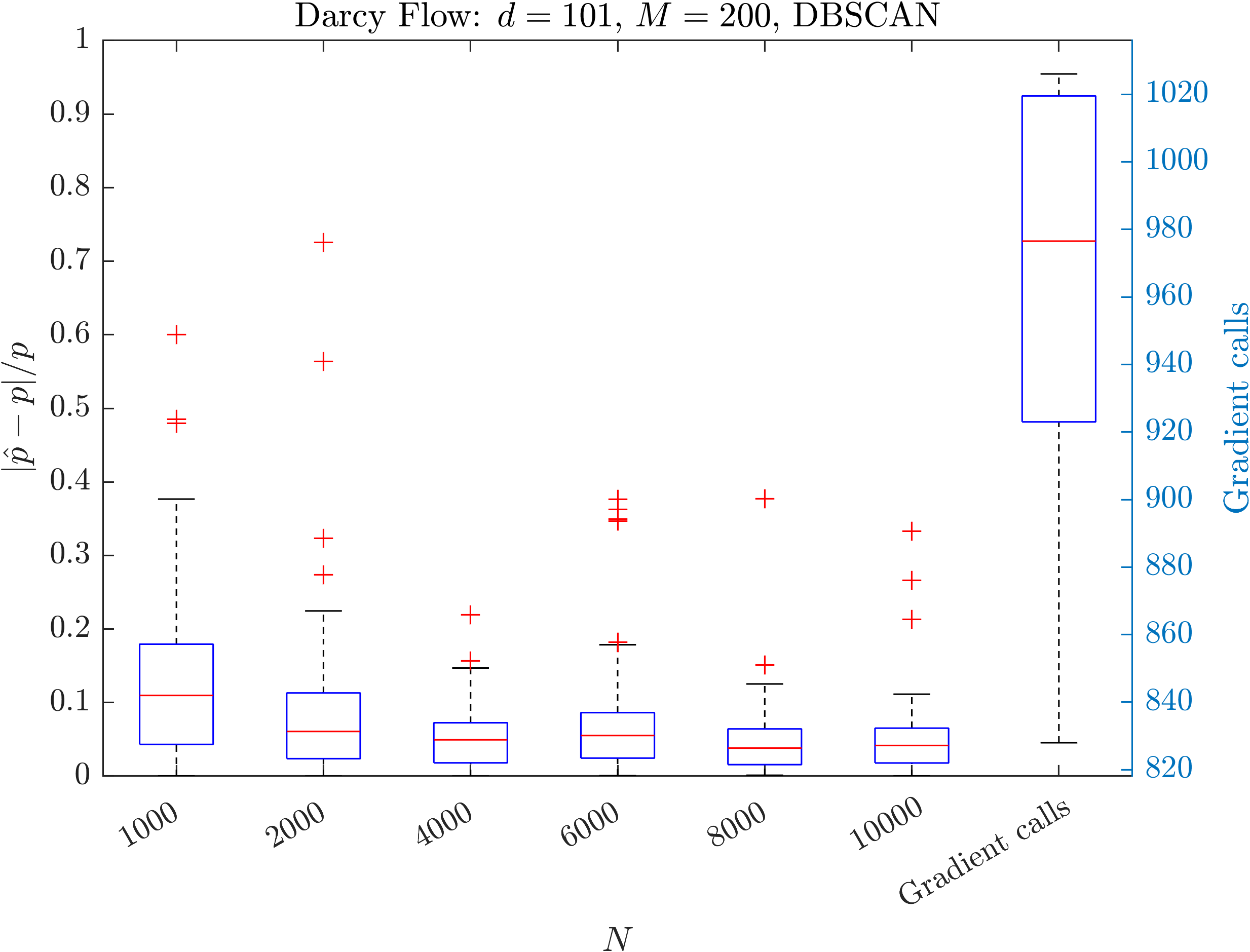}
        \caption{$M=200$} \label{fig:Darcy-2}
    \end{subfigure}
    \caption{Boxplots of normalized absolute error as a function of $N$ and gradient calls of ALDI-IS with DBSCAN for the diffusion problem with $d = 101$, and with respectively $M=100$ and $M=200$ particles in Figure~\ref{fig:Darcy-1} and~\ref{fig:Darcy-2}.}
    \label{fig:Darcy-boxplots}
\end{figure}

In our final experiment, we test ALDI-IS with DBSCAN in the diffusion problem considered in~\cite{thesis_uribe,ehre2024steinvariationalrareevent}, known as the 1D Darcy flow problem. The LSF is defined for any $x \in \R^d$ as
    \begin{align} \label{LSF:Darcy}
    g(x) = 2.7 - \max_{y \in [0,1]}(u(y,x)),
    \end{align}
    where $u(\cdot,x)$ is the solution of 
    \[\frac{\d}{\d y}\left(\kappa(y,x) \frac{\d u}{\d y} (y,x) \right) = - J(y)\,,\quad y \in (0,1),\]
    with the boundary conditions
    \[\frac{\d u}{\d y}(y,x)\biggr\rvert_{y=0} =  - (\sqrt{0.5}x_1 -2), \ u(1,x) = 1\,. \]
    The source term $J$ is defined as
    \[J(y) = 0.8\sum_{i=1}^4 f_{N(0.2i, \ 0.05^2)}(y),\]
    where $f_{N(a, \ b^2)}$ is the PDF of the normal distribution with mean $a$ and standard deviation $b > 0$.
    Lastly, $\kappa$ is modeled as log-normal random field, discretized using a Karhunen-Loève~(KL) expansion with $d - 1$ terms,
    \[\kappa(y,x) = \exp\left(1 + \sqrt{0.3}\sum_{i=1}^{d - 1} \sqrt{\lambda_i} \varphi_i (y) x_{i+1} \right)\]
    where the correlation kernel of the KL expansion is chosen as the exponential kernel such that the eigenpairs $(\lambda_i , \varphi_i (y))$ can be calculated semi-analytically (instead of full numerical resolution by methods in~\cite{betz14}).
    The dimension $d$ of the problem corresponds to the number of terms in the KL expansion + 1 from the random boundary condition on $y=0$. We fix the dimension to $d=101$. In our numerical implementation, the space $[0,1]$ for~$y$ is discretized into 500 equidistant spaced segments. The gradient of the LSF function is approximated by finite difference~(FD), 
    where one approximated gradient call corresponds to two LSF function calls. We therefore regroup the cost of the estimation into total LSF calls, which contains evaluations of $g$ itself plus two times the FD calls. We apply a standard Monte Carlo estimation with $10^9$ samples to obtain a reference probability of $p=7.78\cdot 10^{-6}$. 
    For ALDI with DBSCAN, we use the parameters
    $(q_i) = (1,0.5,0.05,0)$, $(\gamma_i) = (1,0.5,0.01,10^{-3})$, and $(\eps^{\cumu}_i) =(0.1,0.1,0.1,0.005)$.
    The boxplots for the error and LSF gradient call obtained for ALDI-IS with DBSCAN are reported in Figure~\ref{fig:Darcy-boxplots}, in which the scheme shows seemingly good results compared to the LSF gradient calls. Notably, with $M=100$, the performance is quite comparable with $M=200$ but with half the LSF gradient calls.

    We compare the performance with iCE-vMFNM with $1$ initial number of mixture. We regroup the results in Table~\ref{table:Darcy-Flow}, with different numbers of samples per iteration, $N_{\text{it}}$ for iCE-vMFNM and different $N$ and $M$ for ALDI-IS with DBSCAN. We use color codes to highlight the parameters that correspond to comparable total LSF calls between the two methods in order to facilitate comparison. We can see that ALDI-IS with DBSCAN outperforms iCE-vMFNM in our three tested budgets of total LSF calls. However, we mention the drawback of our method, which required more heuristics to tune the parameters, compared to iCE-vMFNM which requires only $N_{\text{it}}$, the initial number of mixture, and the stopping criterion, which is recommended to be set at $0.2$ in most cases (for details, see~\cite{papaioannou_improved_2019}).

\begin{table}
\begin{tabular}{|c|c|c|c|c|c|}
 \hline
Parameters & mean &nRMSE & LSF calls & FD calls & Total LSF calls
 \\ \hline
 \multicolumn{6}{|c|}{iCE-vMFNM}
 \\ \hline $N_{\text{it}} = 500$ & $7.13\cdot 10^{-5}$ & 0.399 &4890 & 0 & 4890
 \\ $N_{\text{it}} = 1000$  & $7.45\cdot 10^{-5}$ &0.213 &8270& 0 &8270
 \\$N_{\text{it}} = 2000$ & $7.41\cdot 10^{-6}$&0.167 &13320& 0 & 13320
 \\ \hline
 \multicolumn{6}{|c|}{ALDI-IS DBSCAN }
 \\ \hline $N=2000, M=200$ &$7.78\cdot 10^{-6}$& 0.261 &2958& 958 & 4874
 \\$N=6000, M=200$ &$7.78\cdot 10^{-6}$ &0.090&6958& 958 &8874
 \\$N=8000, M=400$&$7.89\cdot 10^{-6}$ &0.078&9735& 1735 & 13205
\\ \hline
\end{tabular}
\caption{The nRMSE, mean LSF calls and mean finite difference calls and total LSF calls for iCE-vMFNM and ALDI-IS DBSCAN on the diffusion problem with $d=101$, over $100$ runs.} 
\label{table:Darcy-Flow}
\end{table}
\section{Concluding remarks}
In this paper, we have introduced a new Markov chain importance sampling scheme using ALDI and vMFNM distribution fitting. In Section~\ref{sec:num}, we showed that the scheme coupled with additional numerical strategies is efficient for standard test cases in the literature. The theoretical analysis carried out in this work mainly concerns an idealized and simplified version of ALDI-IS.
Specifically, in Section~\ref{sec:MCIS-theoretical}, we have established the importance sampling error bound that decreases with decreasing $\sigma$ when using i.i.d.\ samples from the intractable $\rhosig$ to estimate $p$. We have also shown that sampling from ULA generates a bias that increases as $\sigma$ decreases. The model fitting error when one fits the PDF of the particles from ULA with a distribution model in order to be employed for importance sampling constitutes part of our ongoing work in the near future.

In Bayesian inference, there exists a gradient-free formulation of the ensemble Kalman filter~\cite{evensen06,Reich_Cotter_2015}, adaptable to the Langevin dynamic~\cite{reich21}. The original formulation is constructed based on the structure of the posterior density and is not straightforwardly applicable to the form of the rare event potential $V_\sigma$. It would be interesting to derive an alternative formulation for $V_\sigma$ to implement a gradient-free algorithm, thereby further saving the gradient calls of the limit state function.

In addition to importance sampling, the Langevin dynamic can be explored in conjunction with subset simulation. 
Instead of introducing the smoothed optimal importance sampling distribution, it would be promising to directly employ a Metropolis adjustment mechanism instead. For example, the well-known Metropolis-adjusted Langevin algorithm (MALA) has been applied to rare event simulation~\cite{pmlr-v206-tit23a}. Similarly, interacting particle sampling methods~\cite{sprungk2025} can be adapted to subset simulation, which will be explored in future work.


\section{Proofs} \label{sec:proofs}
In this section we provide the proofs for Theorem~\ref{thm:error-bound} and Proposition~\ref{prop:all-conditions}.
\subsection{Proof of Theorem~\ref{thm:error-bound}}
    Recall that $X_1 \sim \rhosig$, $Y \sim \rhoopt$ and $Z \sim \rho$. First, we prove this preliminary inequality:
    \begin{align} \label{ineq:largerthanmu}
         &\E\left[F_\sigma(Z) \Indicator{g(Z) > 0}\right] \nonumber 
        \\&= \E\left[\frac{1}{1+\exp\left(\frac{-\mu + g(Z)}{\sigma} \right)} \Indicator{g(Z) \geq 2\mu} \right] + \E\left[\frac{1}{1+\exp\left(\frac{-\mu + g(Z)}{\sigma} \right)} \Indicator{ 0 < g(Z) < 2\mu} \right] \nonumber
        \\&\leq \frac{1}{1+ e^{\mu/\sigma}}+ \E\left[\Indicator{ 0 < g(Z) < 2\mu} \right] \nonumber
        \\&\leq e^{{-\mu}/{\sigma}}+ \P(0 < g(Z) <  2\mu).
    \end{align}
    Since $(X_i)_{i=1,\ldots,N}$ are i.i.d.\ and $\E(\mathds{1}_A(X_i)\rho(X_i) / \rhosig(X_i)) = p$, we have
    \begin{align}\label{main-bound}
        \E\left[\left\lvert \frac{1}{N}\sum_{i=1}^N\frac{ \mathds{1}_A(X_i)\rho(X_i)}{\rhosig(X_i)} - p \right\lvert^2\right]^{1/2} &= \frac{1}{\sqrt{N}} \E\left[\left\lvert \frac{\mathds{1}_A(X_1) \rho(X_1)}{\rhosig(X_1)}- p \right\rvert^2\right]^{1/2} \nonumber
        \\&= \frac{1}{\sqrt{N}} \E\left[\left\lvert \frac{\mathds{1}_A(X_1) p_\sigma}{F_\sigma(X_1)}- p \right\rvert^2\right]^{1/2}\,.
    \end{align}
    Further, the triangular inequality for the $L_2$ norm gives
    \begin{equation} \label{iid-bound}
    \E\left[\left\lvert \frac{\mathds{1}_A(X_1) p_\sigma}{F_\sigma(X_1)}- p \right\rvert^2\right]^{1/2}  \leq \E\left[\left\lvert\frac{\mathds{1}_A(X_1)p_\sigma}{F_\sigma(X_1)} - \frac{p_\sigma}{F_\sigma(Y) }\right\rvert^2\right]^{1/2} + \E\left[\left\lvert\frac{p_\sigma}{F_\sigma(Y) } - p \right\rvert^2\right]^{1/2}.
    \end{equation}
    Note that for any $x \in \R^d$ such that $g(x) \leq 0$, we have$1/F_\sigma(x) = 1 + \exp((-\mu + g(x))/\sigma) \leq 2$ due to $\mu>0$. Then, it holds that $1/F_\sigma(Y) \leq 2$ almost surely, since $g(Y) \leq 0$ for $Y \sim \rhoopt$. For the first term, we can then compute
    \begin{align} \label{bound-firstterm}
     &\E\left[\left\lvert\frac{\mathds{1}_A(X_1)p_\sigma}{F_\sigma(X_1)} - \frac{p_\sigma}{F_\sigma(Y) }\right\rvert^2\right] \nonumber \\& = (p_\sigma)^2 \left(\E\left[\frac{1}{F_\sigma(Y)^2 } \Indicator{ g(X_1) >0}\right] + \E\left[\left\lvert\frac{1}{F_\sigma(X_1)} -\frac{1}{F_\sigma(Y) }\right\rvert^2 \Indicator{ g(X_1) \leq 0}\right]\right) \nonumber
     \\& \leq 4 p_\sigma \P(g(X_1) > 0) + \E\left[\left\lvert e^{(-\mu + g(X_1))/\sigma} - e^{(-\mu + g(Y))/\sigma}\right\rvert^2 \Indicator{ g(X_1) \leq 0}\right] \nonumber
      \\&= 4 p_\sigma \P(g(X_1) > 0) + \E\left[e^{-2\mu/\sigma}\left\lvert e^{g(X_1)/\sigma} -  e^{g(Y)/\sigma}\right\rvert^2 \Indicator{ g(X_1) \leq 0}\right] \nonumber 
     \\&\leq 4 p_\sigma \P(g(X_1) > 0) + \E\left[e^{-2\mu/\sigma}\left\lvert \frac{g(X_1)}{\sigma} -  \frac{g(Y)}{\sigma}\right\rvert^2 \Indicator{ g(X_1) \leq 0}\right] \nonumber 
    \\&= 4 p_\sigma \P(g(X_1) > 0) + \frac{e^{-2\mu/\sigma}}{\sigma^2}\E\left[\left\lvert g(X_1)- g(Y)\right\rvert^2 \Indicator{ g(X_1) \leq 0}\right],
    \end{align}
   where we have used that $\left\lvert e^x - e^y \right\rvert \leq \left\lvert x - y \right\vert$ for $x,y \leq 0$. Remark in addition that (with $Z \sim \rho$), 
    \begin{align*}
        \P(g(X_1) >0) = \E\left[ \frac{\rhosig(Z)}{\rho(Z)}\Indicator{g(Z)> 0}\right] = \frac{1}{p_\sigma} \E\left[ F_\sigma(Z)\Indicator{g(Z) > 0}\right],
    \end{align*}
    so according to~\eqref{ineq:largerthanmu},
    \begin{align*}
        \P(g(X_1) >0) \leq \frac{1}{p_\sigma} \left(e^{{-\mu}/{\sigma}}+ \P(0 < g(Z) <  2\mu)\right),
    \end{align*}
    therefore we get from~\eqref{bound-firstterm} that
       \begin{multline} \label{bound-firstterm-bis}
     \E\left[\left\lvert\frac{\mathds{1}_A(X_1)p_\sigma}{F_\sigma(X_1)} - \frac{p_\sigma}{F_\sigma(Y) }\right \rvert^2 \right] \leq  4 \left(e^{{-\mu}/{\sigma}} + \P(0 < g(Z) <  2\mu) \right)\\
     + \frac{e^{-2\mu/\sigma}}{\sigma^2}\E\left[\left\lvert g(X_1)- g(Y)\right\rvert^2 \Indicator{ g(X_1) \leq 0}\right].
    \end{multline}
    For the second term in the right-hand side of~\eqref{iid-bound}, since $g(Y) \leq 0$ almost surely, we obtain
    \begin{align}\label{bound:second-term}
        \E\left[\left\lvert\frac{p_\sigma}{F_\sigma(Y) } - p \right\rvert^2\right] = \E\left[\left\lvert p_\sigma - p + p_\sigma \exp\left(\frac{-\mu+ g(Y)}{\sigma} \right) \right\rvert^2\right] &\leq 
        \E\left[(\left\lvert p_\sigma - p \right\lvert + \lvert e^{-\mu/\sigma} \rvert)^2\right] \nonumber
        \\&=\left(\left\lvert p_\sigma - p \right\rvert + e^{-\mu/\sigma} \right)^2.
    \end{align}
    Additionally, from~\eqref{eq:smoothed-IS-dist} we have
    \begin{align*}
        &\left\lvert p_\sigma - p \right\rvert = \left\lvert \E\left[F_\sigma(Z)- \mathds{1}_A(Z) \right]\right\rvert 
        \\&\leq \E\left[\left\lvert \frac{1}{1+\exp\left(\frac{-\mu + g(Z)}{\sigma} \right)}- 1 \right\rvert \Indicator{
        g(Z) \leq 0} \right] + \E\left[\frac{1}{1+\exp\left(\frac{-\mu + g(Z)}{\sigma} \right)} \Indicator{ g(Z) > 0 }\right] 
        \\&\leq \E\left[\left\lvert \frac{\exp\left(\frac{-\mu + g(Z)}{\sigma} \right)}{1+\exp\left(\frac{-\mu + g(Z)}{\sigma} \right)} \right\rvert \Indicator{
        g(Z) \leq 0} \right] + e^{{-\mu}/{\sigma}} + \P(0 < g(Z) <  2\mu)\\&\leq 2e^{{-\mu}/{\sigma}} + \P(0 < g(Z) <  2\mu).
    \end{align*}
    where we have applied~\eqref{ineq:largerthanmu}. Injecting this bound into~\eqref{bound:second-term} gives
     \begin{align}\label{bound:second-term-bis}
        \E\left[\left\lvert\frac{p_\sigma}{F_\sigma(Y) } - p \right\rvert^2\right] \leq \left( 3 e^{{-\mu}/{\sigma}} + \P(0 < g(Z) <  2\mu) \right)^2.
    \end{align}
    Lastly, combining the bounds~\eqref{bound-firstterm-bis} and~\eqref{bound:second-term-bis} into~\eqref{iid-bound} and then plugging into~\eqref{main-bound} yields
    \begin{align*}
     &\frac{1}{p}\E\left(\left\lvert \hat p_{\rhosig} - p \right\lvert^2\right)^{1/2} \\
     &\leq \frac{1}{\sqrt{N} p}\biggl[\left( 4 e^{-\mu/\sigma}+ 4 \P(0 < g(Z) <  2\mu) + \frac{e^{-2\mu/\sigma}}{\sigma^2} \E\left[\left\lvert g(X_1)- g(Y)\right\rvert^2 \indicatorA{X_1}\right]\right)^{1/2} 
     \\& \quad \quad \quad \quad + 3 e^{-\mu/\sigma} + \P\left(0 < g(Z) < 2\mu\right)\biggr]
     \\&\leq \frac{3}{\sqrt{N} p}\biggl[\left( e^{-\mu/\sigma}+ \P(0 < g(Z) <  2\mu) + \frac{e^{-2\mu/\sigma}}{\sigma^2}\E\left[\left\lvert g(X_1)- g(Y)\right\rvert^2 \indicatorA{X_1}\right]\right)^{1/2} 
     \\& \quad \quad \quad \quad +e^{-\mu/\sigma} + \P\left(0 < g(Z) < 2\mu\right)\biggr].
\end{align*}
To obtain the bound in the theorem, it remains to note that $a+b \leq 2(a + b + c)^{1/2}$ for $0 \leq a, b \leq 1$ and $c \geq 0$, and to use this inequality with $a = e^{-\mu/\sigma}$, $b = \P(0 < g(Z) < 2 \mu)$ and $c = e^{-2\mu/\sigma} \sigma^{-2} \E[\lvert g(X_1) - g(Y) \rvert^2 \indicatorA{X_1}]^{1/2}$.

\subsection{Proof of Lemma~\ref{lem:sigma-implies-all}}
Assume that~\eqref{cond:sigma-mu} holds and that there exists $\sigma > 0$ such that $\E[g(X_1)^2] < +\infty$.
On one hand this implies $\E[g(Y)^2] < +\infty$ from the inequality
    \begin{align*}
        \int_{\R^d} g^2(x) \rhoopt(x) \, \d x &= \frac{1}{p}\int_{\R^d} g^2(x) \mathds{1}_A(x)\rho(x)\, \d x \\&\leq \frac{1+e^{-\mu/\sigma}}{p}\int_{\R^d} g^2(x) \mathds{1}_A(x) F_\sigma(x) \rho(x)\, \d x
        \\&\leq \frac{1+e^{-\mu/\sigma}}{p}\int_{\R^d} g^2(x) F_\sigma(x) \rho(x) \, \d x
        \\&= \frac{(1+e^{-\mu/\sigma})p_\sigma}{p}\int_{\R^d} g^2(x) \rhosig(x) \, \d x
    \end{align*}
where we have used the inequality $1 \leq F_\sigma(x)(1+e^{-\mu/\sigma})$ for any $x \in \R^d$ such that $g(x) \leq 0$.  On the other hand, we have also $\E[g(X_1)^2] < +\infty$ for any $\sigma > 0$, since for any $\sigma_1, \sigma_2 > 0$,
    \begin{align*}
        \int_{\R^d} g^2(x) \rhoopt_{\sigma_1}(x) \, \d x &= \frac{1}{p_{\sigma_1}}\int_{\R^d} g^2(x) F_{\sigma_1}(x)\rho(x)\, \d x 
        \\&\leq \frac{1+e^{-\mu/{\sigma_2}}}{p_{\sigma_1}}\int_{\R^d} g^2(x)  F_{\sigma_1}(x)F_{\sigma_2}(x) \rho(x)\, \d x
        \\&\leq \frac{1+e^{-\mu/{\sigma_2}}}{p_{\sigma_1}}\int_{\R^d} g^2(x) F_{\sigma_2}(x) \rho(x) \, \d x
        \\&= \frac{(1+e^{-\mu/{\sigma_2}})p_{\sigma_2}}{p_{\sigma_1}}\int_{\R^d} g^2(x) \rhoopt_{\sigma_2}(x) \, \d x.
    \end{align*}
    At last, we also have $\lim_{\sigma \to 0} \E[g(X_1)^2] = \E[g(Y)^2] < +\infty$ by Lebesgue's dominated convergence theorem, since for any $x \in \R^d$, $\lim_{\sigma \to 0}g(x)^2 \rhosig(x) \to g(x)^2 \rhoopt(x) $ and 
    \begin{align*}
        g(x)^2 \rhosig(x) = g(x)^2\frac{F_\sigma(x)\rho(x)}{p_\sigma(x)}    \leq \frac{2}{p} g(x)^2 \rho(x)
    \end{align*}
    which is integrable. In the previous display we have used the fact that
    \[p_\sigma = \int F_\sigma(x)\rho(x) \d x \geq \int \Indicator{g(x) \leq 0}F_\sigma(x) \rho (x) \d x \geq \int \Indicator{g(x) \leq 0}\frac{1}{2} \rho (x) \d x.\]
\subsection{Proof of Proposition~\ref{prop:all-conditions}}
To verify~\eqref{cond:log-cc-log-smooth}, the gradient and Hessian of $V_\sigma$ are required and are respectively given for any $x \in \R^d$ by  
\begin{align}\label{eq:grad-V}
    \grad V_\sigma(x) = -\nabla \log F_\sigma(x) + x  = \frac{\nabla g(x)}{\sigma}\left(1- F_\sigma(x)\right) + x\,,    
\end{align}
and
\begin{align}\label{eq:Hess-V}
\Hess(V_\sigma)(x) &= \Hess(-\log F_\sigma)(x) + I \nonumber \\&= \frac{(1-F_\sigma(x))F_\sigma(x)}{\sigma^2}\nabla g(x) \nabla g(x)^\top + \frac{1-F_\sigma(x)}{\sigma}\Hess (g)(x) + I.
\end{align}

First, we verify that $H_\sigma^g \leq G$. Since $\nabla g$ is $G$-Lipschitz continuous, we have $\Hess(g)(x) \leq GI$ for any $x \in \R^d$. We also have $1-F_\sigma(x) \leq 1$ for any $x \in \R^d$ from Definition~\ref{def:smooth-indicator}. So $H_\sigma^g = \inf_{x \in \R^d}(1-F_\sigma(x))\lambda_{\min}(\Hess(g)(x)) \leq G$.

Let $x \in \R^d$. The function $g$ is of class $C^2(\R^d,\R)$ so from Schwarz's theorem, $\Hess(-\log F_\sigma)(x)$ is symmetric and according to~\eqref{eq:Hess-V},
\[\lambda_{\min}(\Hess (V_\sigma)(x)) = \lambda_{\min}(\Hess (- \log F_\sigma)(x)) + 1.\]
Note that the only eigenvalues of $\nabla g(x) \nabla g(x)^\top$ are $\lVert \grad g(x)\rVert^2$ and $0$, so $\lambda_{\min}(\nabla g(x) \nabla g(x)^\top) = 0$. By Weyl's inequality on~\eqref{eq:Hess-V},
\[\lambda_{\min}(\Hess (V_\sigma)(x)) \geq 0 + \frac{1 - F_\sigma(x)}{\sigma}\lambda_{\min}(\Hess (g)(x)) + 1.\]
Observe then that the definition of $\alpha_\sigma$ gives $\inf_{x \in \R^d}\lambda_{\min}(\Hess (V_\sigma)(x)) \geq \alpha_\sigma$, so for any $x \in \R^d, \Hess(V_\sigma)(x) \geq \alpha_\sigma I$.
To prove that $\Hess(V_\sigma)(x) \leq L_\sigma I$ for any $x \in \R^d$, it is equivalent to prove that for any $x,y \in \R^d$,
\[ \lVert \nabla V_\sigma(x) - \nabla V_\sigma(y) \rVert \leq L_\sigma \lVert x - y \rVert.\]
Let $x, y \in \R^d$,
\begin{align*}
    &\lVert \nabla V_\sigma(x) - \nabla V_\sigma(y) \rVert \\&=  \lVert - \nabla \log F_\sigma(x) + \nabla \log F_\sigma(y) + ( x - y) \rVert 
    \\&\leq \left\lVert \frac{1}{\sigma}(\nabla g(x) - \nabla g(y)) - F_\sigma(x)\frac{\nabla g(x)}{\sigma} + F_\sigma(y)\frac{\nabla g(y)}{\sigma} \right\rVert + \lVert x - y \rVert
    \\&\leq \frac{1}{\sigma}\left(\left\lVert \nabla g(x) - \nabla g(y)\right\rVert + \left\lVert F_\sigma(x)\nabla g(x) - F_\sigma(y)\nabla g(y)\right\rVert\right)  + \lVert x - y \rVert.
\end{align*}
In addition,
\begin{align*}
    &\left\lVert F_\sigma(x)\nabla g(x) - F_\sigma(y)\nabla g(y)\right\rVert \\&= \left\lVert F_\sigma(x)\nabla g(x) -F_\sigma(x)\nabla g(y) + F_\sigma(x) \nabla g(y) - F_\sigma(y)\nabla g(y)\right\rVert
    \\&\leq \lvert F_\sigma(x)\rvert\left\lVert \nabla g(x) - \nabla g(y)\right\rVert + \lVert \nabla g(y) \rVert \left \lvert F_\sigma(x) - F_\sigma(y) \right \rvert 
    \\&\leq G\left\lVert x - y\right\rVert + K \left \lvert F_\sigma(x) - F_\sigma(y) \right \rvert 
\end{align*}
where we have used the facts that $g$ is $G$-smooth and $K$-Lipschitz continuous, the latter implying that $\lVert \grad g \rVert$ is bounded by $K$. We also have 
\begin{align*}
    \lvert F_\sigma(x) - F_\sigma(y) \rvert
    &= \left \lvert \left(1 + e^{\left(-\mu + g(x)\right)/\sigma} \right)^{-1} - \left(1 + e^{(-\mu + g(y))/\sigma} \right)^{-1}\right \rvert
    \\&\leq \left \lvert e^{\left(-\mu + g(x)\right)/\sigma}  - e^{(-\mu + g(y))/\sigma} \right \rvert
    \\&= e^{-\mu/\sigma}\left \lvert e^{g(x)/\sigma}  - e^{g(y)/\sigma} \right \rvert
    \\&\leq e^{-\mu/\sigma}e^{r/\sigma}\left \lvert g(x)/\sigma  - g(y)/\sigma \right \rvert
    \\&\leq \frac{e^{(-\mu + r)/\sigma}K}{\sigma}\left \lVert x - y \right \rVert
\end{align*}
where we have used the facts that the exponential function is $e^{r/\sigma}$-Lipschitz continuous on $[-r/\sigma, r/\sigma]$ and $g$ is $K$-Lipschitz continuous. Collecting all the inequalities yields 
\begin{align*}
    &\lVert \nabla V_\sigma(x) - \nabla V_\sigma(y) \rVert 
    \\&\leq \frac{1}{\sigma}\left(\left\lVert \nabla g(x) - \nabla g(y)\right\rVert + \left\lVert F_\sigma(x)\nabla g(x) - F_\sigma(y)\nabla g(y)\right\rVert\right)  + \lVert x - y \rVert
    \\&\leq \left(\frac{2G}{\sigma} + \frac{e^{(-\mu + r)/\sigma}}{\sigma^2}K^2+1\right)\lVert x - y \rVert.
\end{align*}

\section*{Acknowledgements}
The author J. Beh is enrolled in a Ph.D. program co-funded by
ONERA - The French Aerospace Lab and the University Research School EUR-MINT (State
support managed by the National Research Agency for Future Investments program bearing the reference ANR-18-EURE-0023). Their financial supports are gratefully acknowledged. J. Beh would also like to express their gratitude towards the Mathematisches Institut of Universität Mannheim for their stay in Mannheim, as well as \'Ecole Doctorale Mathématiques, Informatique, Télécommunications de Toulouse and Université Toulouse III - Paul Sabatier for financially supporting the stay.
\printbibliography
\end{document}